\documentclass[11pt]{article}

\usepackage{amssymb}
\usepackage{MnSymbol}
\newfont{\boldit}{cmbxti10}
\def\La{\Lambda}
\def\In{I}  

\def\ell{{{l}}}
\def\dimk{{\rm dim}_k}

\def\la{\lambda}

\def\mod{{\rm mod}}
\def\odimk{\underline{\rm dim}\,}

\def\crk{{\rm cor}}
\def\r{{\rm r}}
\def\Hom{{\rm Hom}}
\def\trivm{\varepsilon}
\def\ma21#1#2{{\tiny\left[\begin{array}{c}\!\!\!\!#1\!\!\!\!\!\\\!\!\!#2\!\!\!\!\end{array}\right]}} 
\def\wmm#1#2{{\footnotesize\left[\begin{array}{c}\!\!\!\!#1\!\!\!\!\!\\\!\!\!#2\!\!\!\!\\\end{array}\right]}} 
\def\wmmm#1#2#3{{\footnotesize\left[\begin{array}{c}\!\!\!\!#1\!\!\!\!\!\\\!\!\!#2\!\!\!\!\\\!\!\!#3\!\!\!\!\end{array}\right]}} 
\def\wmmmm#1#2#3#4{{\footnotesize\left[\begin{array}{c}\!\!\!\!#1\!\!\!\!\!\\\!\!\!#2\!\!\!\!\\\!\!\!#3\!\!\!\!\\\!\!\!#4\!\!\!\!\end{array}\right]}} 

\def\cor{{\rm cor}}

\def\sk{\smallskip} 

\def\ZRlla{Z_{R(l,\la)}}
\def\ZIm0{Z_{I(m,0)}}
\def\ZpIm0{Z'_{I(m,0)}}
\def\ZPm0{Z_{P(m,0)}}
\def\ZpPm0{Z'_{P(m,0)}}
\def\ZP2m11{Z_{P(2m+1,1)}}
\def\ZPdmjd{Z_{P(2m+1,2)}}
\def\ZPdmjt{Z_{P(2m+1,3)}}
\def\ZPdmjc{Z_{P(2m+1,4)}}
\def\ZRzdlz{Z_{R(0,2l,0)}}
\def\ZRjdlz{Z_{R(1,2l,0)}}
\def\ZRzdlj{Z_{R(0,2l,1)}}
\def\ZRjdlj{Z_{R(1,2l,1)}}
\def\ZRzdli{Z_{R(0,2l,\infty)}}
\def\ZRjdli{Z_{R(1,2l,\infty)}}
\def\ZPdmj{Z_{P(2m,1)}}
\def\ZpPdmj{Z'_{P(2m,1)}}
\def\ZPdmd{Z_{P(2m,2)}}
\def\ZpPdmd{Z'_{P(2m,2)}}
\def\ZPdmt{Z_{P(2m,3)}}
\def\ZpPdmt{Z'_{P(2m,3)}}
\def\ZPdmc{Z_{P(2m,4)}}
\def\ZpPdmc{Z'_{P(2m,4)}}
\def\ZIdmj{Z_{I(2m,1)}}
\def\ZpIdmj{Z'_{I(2m,1)}}
\def\ZIdmd{Z_{I(2m,2)}}
\def\ZpIdmd{Z'_{I(2m,2)}}
\def\ZIdmt{Z_{I(2m,3)}}
\def\ZpIdmt{Z'_{I(2m,3)}}
\def\ZIdmc{Z_{I(2m,4)}}
\def\ZpIdmc{Z'_{I(2m,4)}}
\def\ZIdmjj{Z_{I(2m+1,1)}}
\def\ZIdmjd{Z_{I(2m+1,2)}}
\def\ZIdmjt{Z_{I(2m+1,3)}}
\def\ZIdmjc{Z_{I(2m+1,4)}}
\def\ZRzdljj{Z_{R(0,2l-1,1)}}
\def\ZpRzdljj{Z'_{R(0,2l-1,1)}}
\def\ZRjdljj{Z_{R(1,2l-1,1)}}
\def\ZpRjdljj{Z'_{R(1,2l-1,1)}}
\def\ZRzdljz{Z_{R(0,2l-1,0)}}
\def\ZpRzdljz{Z'_{R(0,2l-1,0)}}
\def\ZRjdljz{Z_{R(1,2l-1,0)}}
\def\ZpRjdljz{Z'_{R(1,2l-1,0)}}
\def\ZRzdlji{Z_{R(0,2l-1,\infty)}}
\def\ZpRzdlji{Z'_{R(0,2l-1,\infty)}}
\def\ZRjdlji{Z_{R(1,2l-1,\infty)}}
\def\ZpRjdlji{Z'_{R(1,2l-1,\infty)}}
\def\ZRzdlj{Z_{R(0,2l,1)}}
\def\ZRjdlj{Z_{R(1,2l,1)}}
\def\ZRzdlz{Z_{R(0,2l,0)}}
\def\ZRjdlz{Z_{R(1,2l,0)}}
\def\ZRzdli{Z_{R(0,2l,\infty)}}
\def\ZRjdli{Z_{R(1,2l,\infty)}}

\def\DM{\CM_2} 
\def\DML{\CM_1} 
\def\DN{\CM_3}

\def\wlj{{^\circ\!\pi}} 
\def\wld{\pi^\circ}
\def\Iod{\overline{I}}

\def\MM{\mathbb M}
\def\M{\mathbb M}

\def\PP{\mathbb P}

\def\ZZ{\mathbb Z}

\def\CT{{\cal  T}}
\def\CM{{\cal M}}
\def\CN{{\cal N}}

\def\CX{{\cal X}}

\makeatletter
\renewcommand\normalsize{%
   \@setfontsize\normalsize\@xipt{13.6}%
   \abovedisplayskip 11\p@ \@plus3\p@ \@minus6\p@
   \abovedisplayshortskip \z@ \@plus3\p@
   \belowdisplayshortskip 6.5\p@ \@plus3.5\p@ \@minus3\p@
   \belowdisplayskip \abovedisplayskip
   \let\@listi\@listN}

\renewcommand\small{%
   \@setfontsize\small\@xpt\@xiipt
   \abovedisplayskip 10\p@ \@plus2\p@ \@minus5\p@
   \abovedisplayshortskip \z@ \@plus3\p@
   \belowdisplayshortskip 6\p@ \@plus3\p@ \@minus3\p@
   \belowdisplayskip \abovedisplayskip
   \let\@listi\@listS}

\renewcommand\footnotesize{%
   \@setfontsize\footnotesize\@ixpt{11}%
   \abovedisplayskip 8\p@ \@plus2\p@ \@minus4\p@
   \abovedisplayshortskip \z@ \@plus\p@
   \belowdisplayshortskip 4\p@ \@plus2\p@ \@minus2\p@
   \belowdisplayskip \abovedisplayskip
   \let\@listi\@listF}
\makeatother

\makeatletter
\def\@listN{%
        \listparindent 1.5em%
        \leftmargini 3.6em%
        \leftmarginii 1.5em%
        \leftmarginiii 1.5em%
        \leftmarginiv 1.5em%
        \leftmarginv 1.5em%
        \leftmarginvi 1.5em%
        \labelsep .5em%
        \labelwidth\leftmargini%
        \advance\labelwidth-\labelsep%
        \leftmargin\leftmargini
        \topsep 3.5pt plus 0.5pt minus 0.5pt%
        \parsep 2.0pt plus 0.2pt minus 0.1pt%
        \partopsep 0.5pt plus 0.5pt minus 0.1pt%
        \itemsep 2.0pt plus 0.2pt minus 0.1pt}%
\def\@listS{%
        \listparindent 1.5em%
        \leftmargini 3.6em%
        \leftmarginii 1.5em%
        \leftmarginiii 1.5em%
        \leftmarginiv 1.5em%
        \leftmarginv 1.5em%
        \leftmarginvi 1.5em%
        \labelsep .5em%
        \labelwidth\leftmargini%
        \advance\labelwidth-\labelsep%
        \leftmargin\leftmargini%
        \topsep 3pt plus 0.2pt minus 0.2pt%
        \parsep 2.0pt plus 0.2pt minus 0.1pt%
        \partopsep 0.5pt plus 0.2pt minus 0.1pt%
        \itemsep \parsep}%
\def\@listF{%
        \listparindent 1.5em%
        \leftmargini 3.6em%
        \leftmarginii 1.5em%
        \leftmarginiii 1.5em%
        \leftmarginiv 1.5em%
        \leftmarginv 1.5em%
        \leftmarginvi 1.5em%
        \labelsep .5em%
        \labelwidth\leftmargini%
        \advance\labelwidth-\labelsep%
        \leftmargin\leftmargini%
        \topsep 3pt plus 0.0pt minus 0.0pt%
        \parsep 2.0pt plus 2.0pt minus 0.0pt%
        \partopsep 0.0pt plus 0.0pt minus 0.0pt%
        \itemsep \parsep}%
\def\@listii{%
        \listparindent 1em%
        \leftmargin\leftmarginii%
        \labelwidth\leftmarginii%
        \labelsep .5em%
        \advance\labelwidth-\labelsep%
        \topsep 1pt plus 0.1pt minus 0pt%
        \parsep 0.0pt plus 0.1pt minus 0pt%
        \partopsep 0.5pt
        \itemsep \parsep}%
\def\@listiii{%
        \listparindent 1em%
        \leftmargin\leftmarginiii%
        \labelwidth\leftmarginiii%
        \labelsep .5em%
        \advance\labelwidth-\labelsep%
        \topsep 0pt
        \parsep 0pt
        \partopsep 0pt
        \itemsep \topsep}%
\def\@listiv{%
        \listparindent 1em%
        \leftmargin\leftmarginiv%
        \labelwidth\leftmarginiv%
        \labelsep .5em%
        \advance\labelwidth-\labelsep%
        \topsep 0pt%
        \parsep 0pt%
        \partopsep 0pt%
        \itemsep \topsep}%
\def\@listv{%
        \listparindent 1em%
        \leftmargin\leftmarginv%
        \labelwidth\leftmarginv%
        \labelsep .5em%
        \advance\labelwidth-\labelsep%
        \topsep 0pt%
        \parsep 0pt%
        \partopsep 0pt%
        \itemsep \topsep}%
\def\@listvi{%
        \listparindent 1em%
        \leftmargin\leftmarginvi%
        \labelwidth\leftmarginvi%
        \labelsep .5em%
        \advance\labelwidth-\labelsep%
        \topsep 0pt%
        \parsep 0pt%
        \partopsep 0pt%
        \itemsep \topsep}%
\normalsize \makeatother

\begin{document}

\begin{center}
{\Large\sc The dimensions of the homomorphism spaces to
indecomposable modules over the four subspace algebra
}\\

\medskip
{\sc by}\\
\medskip
{ \sc  Andrzej Mr\'oz (TORU\'N, 2010)}

 {}\footnotetext[1]{AMS
 2000
 Classification: 16G60, 16G70, 14L30, 68Q99.}
\footnotetext[2]{Key words: four subspace algebra, indecomposable
module, homomorphism space isomorphism problem, multiplicity
problem.}

\end{center}

\vspace{2mm}

\bigskip

\indent {\footnotesize \bf Abstract.}
 {\footnotesize
 This is the addendum to the paper \cite{df}. We give here the
 full proof of \cite[Proposition 3.3]{df}, describing the formulas for the dimensions of the homomorphism spaces to indecomposable modules over the four subspace algebra $\La$.
  The reader can also find here the full description of
 indecomposable $\La$-modules.
 }

\bigskip

{\bf Introduction and preliminaries.} We briefly recall the
setting from \cite{df}; for more details we refer to that article
and also to \cite{ASS, GePo2, Siks, SS}.

By $k$ we denote the fixed algebraically closed field and by
$\MM_{n\times m}$, the set of all $(n\times m)$-matrices over $k$,
for $n,m\geq 0$ (we also consider the trivial matrices with zero
rows or columns). Given  $A\in \M_{m\times n}$, we denote by ${\rm
r}(A)$ the rank of $A$ and by $\crk(A)=m-{\r}(A)$ its
 corank.

\smallskip

 Let $Q$ be
the quiver
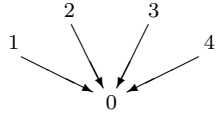
\begin{center}
\begin{picture}(50,40)
\put(-7,20){\vector(2,-1){27}} \put(-12,23){\scriptsize $1$}

\put(9,35){\scriptsize $2$} \put(12,32){\vector(1,-2){12}}

\put(41,35){\scriptsize $3$} \put(41,32){\vector(-1,-2){12}}

\put(62,23){\scriptsize $4$} \put(60,20){\vector(-2,-1){27}}

\put(25,0){\scriptsize $0$}
\end{picture}
\end{center}
(the set of vertices $\{0,1,2,3,4\}$ we denote by $Q_0$). The path
algebra $\La=kQ$ we call the {\em four subspace algebra}. Then any
collection of matrices $(A,B,C,D)\in\MM_{n_0\times n_1}\times
\MM_{n_0\times n_2}\times\MM_{n_0\times n_3}\times\MM_{n_0\times
n_4}$, $n_i\geq 0$, $i=0,\ldots,4$, determines $\La$-module
(treated as the matrix representation of the quiver $Q$) of the
form
\smallskip

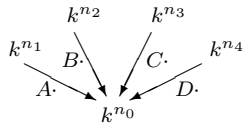
\begin{center}
\begin{picture}(50,40)
\put(-7,20){\vector(2,-1){27}} \put(-13,23){\scriptsize $k^{n_1}$}
\put(-3,8){\scriptsize $A\cdot$}

\put(9,36){\scriptsize $k^{n_2}$} \put(12,32){\vector(1,-2){12}}
\put(7,19){\scriptsize $B\cdot$}

\put(41,36){\scriptsize $k^{n_3}$} \put(41,32){\vector(-1,-2){12}}
\put(39,19){\scriptsize $C\cdot$}

\put(63,23){\scriptsize $k^{n_4}$} \put(60,20){\vector(-2,-1){27}}
\put(50,8){\scriptsize $D\cdot$}

\put(22,-2){\scriptsize $k^{n_0}$}
\end{picture}
\end{center}
 Further on we will present $\La$-modules as the
quadruples $(A,B,C,D)$ as above. For $\La$-module $M$, the
dimension vector $\odimk(M)$ we will present in the form
$$[\dimk(M_0),\dimk(M_1),\dimk(M_2),\dimk(M_3),\dimk(M_4)].$$
By $\tau=\tau_\La$ we denote the Auslander-Reiten translate in the
category $\mod\,\La$ of finite-dimensional $\La$-modules.

\smallskip

 We
recall that all the indecomposable $\La$-modules  are divided into
the following three families:
\begin{itemize}
\item postprojective modules of the form $P(m,j):=\tau^{-m}P(j)$,
for all $m\geq 0$, $j\in Q_0$, where $P(i)$ is the indecomposable
projective $\La$-module corresponding to the vertex $i$;
\item preinjective modules of the form $\In(m,j):=\tau^m \In(j)$,
for all $m\geq 0$, $j\in Q_0$, where $I(i)$ is the indecomposable
injective $\La$-module corresponding to the vertex $i$;
\item regular  modules forming a family of pairwise orthogonal stable standard
tubes $\CT_\la$ in Auslaned-Reiten quiver of the category
$\mod\,\La$, for all $\la\in\PP^1(k)$, where
$\CT_0,\CT_1,\CT_\infty$ are the tubes of rank 2 and the remaining
are homogeneous. The modules from the tubes $\CT_\la$, for
$\la\in\{0,1,\infty\}$, are denoted by $R(s, l,\la)$, where
$s\in\ZZ_2$ is the number of quasi-socle, $l\geq 1$ is
quasi-length (that is,  $R(s,l,\la)$ has quasi-length $l$, its
quasi-socle is isomorphic to $R(s,1,\la)$, and moreover $\tau
R(s,l,\la) \cong R(s\oplus_2 1,l,\la)$). The modules from
homogeneous tubes $\CT_\la$, for  $\la\in k\setminus\{0,1\}$, we
denote by $R(l,\la)$, where $l\geq 1$ is quasi-length (in
particular, $\tau R(l,\la)=R(l,\la)$).

\end{itemize}

\medskip

The paper \cite{df} deals with two natural problems: the {\em
multiplicity problem} (the algorithm determining the
 multiplicities of all the indecomposables appearing in the
 indecomposable decomposition of $\La$-module) and the {\em isomorphism problem} (the algorithmic
 criterion for deciding if two $\La$-modules are isomorphic).

 In
 the above
 considerations the crucial role is played by the direct
 formulas for computing the dimensions
 $[M,X]:=\dimk(\Hom_\La(M,X))$ of homomorphism spaces from given
 $\La$-module $M$ to arbitrary indecomposable $\La$-module $X$.
 Such formulas were given in \cite[Proposition 3.3]{df}.
 The aim of this paper is to give the full proof of the
 proposition, which we recall below in 1.

In the proof, given in 3, we use the matrix description of the
indecomposable $\La$-modules, which we recall in Proposition 2.

\bigskip

{\bf 1.} We recall below the assertion of \cite[Proposition
3.3]{df}. First we establish some notation. For the integers
$x,y,x',y',u,v,s,t,i \geq 0$ with $u\leq x$, $v\leq y,y'$, $s\leq
x,x'$, $t\leq y$ and the matrices $Z\in\MM_{x\times y}$,
$Z'\in\MM_{x'\times y'}$, $F\in\MM_{u\times v}$, $E\in\MM_{s\times
t}$ we define the block matrices $\DML^i=\DML(Z,Z',F,i),
\DM^i=\DM(Z,Z',E,i)\in \MM_{(x'+ix)\times (y'+iy)}$ and
$\DN^i=\DN(Z,Z',E,i)\in\MM_{(x'+ix)\times (y'+(i+1)y)}$ as
follows:
$$\DML^i:=\left[{\scriptsize\begin{array}{ccccc} Z'&&&&\\
F''&Z&&\\
&F'&Z&&\\
&&\ddots&\ddots&\\
&&&F'&Z
\end{array}}\right],\ \ \ \ \ \DM^i:=\left[{\scriptsize\begin{array}{ccccc} Z'&E''&&&\\
&Z&E'&\\
&&\ddots&\ddots&\\
&&&Z&E'\\
&&&&Z
\end{array}}\right],$$
$$\DN^i:=\left[{\scriptsize\begin{array}{cccccc} Z'&E''&&&&\\
&Z&E'&&\\
&&\ddots&\ddots&&\\
&&&Z&E'&\\
&&&&Z&E'
\end{array}}\right]$$
(the blocks $Z'$, $E''$, $F''$ appear only once, the remaining
coefficients are
zero),  where $F'=\left[{\scriptsize\begin{array}{cc} 0&F\\
0&0\end{array}}\right], E'=\left[{\scriptsize\begin{array}{cc} 0&0\\
E&0\end{array}}\right] \in\MM_{x\times y}$, $F''=\left[{\scriptsize\begin{array}{cc} 0&F\\
0&0\end{array}}\right]\in\MM_{x\times y'}$ and $E''=\left[{\scriptsize\begin{array}{cc} 0&0\\
E&0\end{array}}\right]\in\MM_{x'\times y}$. In particular,
$\DML^0=\DM^0=Z'$ and $\DN^0=[Z'|E'']$.

For $A\in\MM_{x\times y}$, $B\in\MM_{z\times t}$, by $A\oplus B$
we mean the matrix ${\scriptsize
\left[\begin{array}{cc}A&\!0\\0&\!B
\end{array}\right]}\in\MM_{(x+z)\times (y+t)}$.

\bigskip
\newpage
\def\ZRlla{Z_{R(l,\la)}}
\def\ZIm0{Z_{I(n,0)}}
\def\ZpIm0{Z'_{I(n,0)}}
\def\ZPm0{Z_{P(n,0)}}
\def\ZpPm0{Z'_{P(n,0)}}
\def\ZP2m11{Z_{P(2n+1,1)}}
\def\ZPdmjd{Z_{P(2n+1,2)}}
\def\ZPdmjt{Z_{P(2n+1,3)}}
\def\ZPdmjc{Z_{P(2n+1,4)}}
\def\ZRzdlz{Z_{R(0,2l,0)}}
\def\ZRjdlz{Z_{R(1,2l,0)}}
\def\ZRzdlj{Z_{R(0,2l,1)}}
\def\ZRjdlj{Z_{R(1,2l,1)}}
\def\ZRzdli{Z_{R(0,2l,\infty)}}
\def\ZRjdli{Z_{R(1,2l,\infty)}}
\def\ZPdmj{Z_{P(2n,1)}}
\def\ZpPdmj{Z'_{P(2n,1)}}
\def\ZPdmd{Z_{P(2n,2)}}
\def\ZpPdmd{Z'_{P(2n,2)}}
\def\ZPdmt{Z_{P(2n,3)}}
\def\ZpPdmt{Z'_{P(2n,3)}}
\def\ZPdmc{Z_{P(2n,4)}}
\def\ZpPdmc{Z'_{P(2n,4)}}
\def\ZIdmj{Z_{I(2n,1)}}
\def\ZpIdmj{Z'_{I(2n,1)}}
\def\ZIdmd{Z_{I(2n,2)}}
\def\ZpIdmd{Z'_{I(2n,2)}}
\def\ZIdmt{Z_{I(2n,3)}}
\def\ZpIdmt{Z'_{I(2n,3)}}
\def\ZIdmc{Z_{I(2n,4)}}
\def\ZpIdmc{Z'_{I(2n,4)}}
\def\ZIdmjj{Z_{I(2n+1,1)}}
\def\ZIdmjd{Z_{I(2n+1,2)}}
\def\ZIdmjt{Z_{I(2n+1,3)}}
\def\ZIdmjc{Z_{I(2n+1,4)}}
\def\ZRzdljj{Z_{R(0,2l-1,1)}}
\def\ZpRzdljj{Z'_{R(0,2l-1,1)}}
\def\ZRjdljj{Z_{R(1,2l-1,1)}}
\def\ZpRjdljj{Z'_{R(1,2l-1,1)}}
\def\ZRzdljz{Z_{R(0,2l-1,0)}}
\def\ZpRzdljz{Z'_{R(0,2l-1,0)}}
\def\ZRjdljz{Z_{R(1,2l-1,0)}}
\def\ZpRjdljz{Z'_{R(1,2l-1,0)}}
\def\ZRzdlji{Z_{R(0,2l-1,\infty)}}
\def\ZpRzdlji{Z'_{R(0,2l-1,\infty)}}
\def\ZRjdlji{Z_{R(1,2l-1,\infty)}}
\def\ZpRjdlji{Z'_{R(1,2l-1,\infty)}}
\def\ZRzdlj{Z_{R(0,2l,1)}}
\def\ZRjdlj{Z_{R(1,2l,1)}}
\def\ZRzdlz{Z_{R(0,2l,0)}}
\def\ZRjdlz{Z_{R(1,2l,0)}}
\def\ZRzdli{Z_{R(0,2l,\infty)}}
\def\ZRjdli{Z_{R(1,2l,\infty)}}

{\bf Proposition.} {\em Let $M=(A,B,C,D)$ be a fixed $\La$-module,
with $A\in\MM_{n_0\times n_1}$, $B\in\MM_{n_0\times n_2}$,
$C\in\MM_{n_0\times n_3}$, $D\in\MM_{n_0\times n_4}$. Then the
formulae for the dimension $[M,X]$, for indecomposable
$\La$-module $X\in\CX$, are given in the following table
\begin{center}
\begin{tabular}{ll}
\hline $X$&$[M,X]$\smallskip\\
\hline\hline $P(0,0)$&$\crk([A\ B\ C\ D])$\sk\\
\hline
$P(n,0)$,& $\crk(\DN(\ZPm0,\ZpPm0,C\oplus D,n-1)$,\sk\\
$n\geq 1$&$\ZpPm0=\left[{\scriptsize\begin{array}{cccccc}A&\!\!\!0&\!\!\!B&\!\!\!0&\!\!\!C&\!\!\!D\\0&\!\!\!0&\!\!\!0&\!\!\!B&\!\!\!0&\!\!\!-D\\
0&\!\!\!A&\!\!\!0&\!\!\!0&\!\!\!-C&\!\!\!0\\\end{array}}\right]$, $\ZPm0=\left[{\scriptsize\begin{array}{cccc}0&\!\!\!-D&\!\!\!0&\!\!\!B\\-C&\!\!\!0&\!\!\!A&\!\!\!0\end{array}}\right]$\sk\\
\hline\hline $P(2n+1,1)$,& $\crk(\DN(\ZP2m11,\ZP2m11,-A,n)),$\sk\\
 $n\geq 0$&  $\ZP2m11=\left[{\scriptsize\begin{array}{cccc}A&\!\!\!0&\!\!\!C&\!\!\!D\\0&\!\!\!B&\!\!\!0&\!\!\!-D\end{array}}\right]$\sk\\
 \hline $P(2n+1,2)$,& $\crk(\DN(\ZPdmjd,\ZPdmjd,-B,n)),$\sk\\
 $n\geq 0$&  $\ZPdmjd=\left[{\scriptsize\begin{array}{cccc}B&\!\!\!0&\!\!\!D&\!\!\!A\\0&\!\!\!C&\!\!\!0&\!\!\!-A\end{array}}\right]$\sk\\
 \hline $P(2n+1,3)$,& $\crk(\DN(\ZPdmjt,\ZPdmjt,-C,n)),$\sk\\
 $n\geq 0$&  $\ZPdmjt=\left[{\scriptsize\begin{array}{cccc}C&\!\!\!0&\!\!\!A&\!\!\!B\\0&\!\!\!D&\!\!\!0&\!\!\!-B\end{array}}\right]$\sk\\
 \hline $P(2n+1,4)$,& $\crk(\DN(\ZPdmjc,\ZPdmjc,-D,n)),$\sk\\
 $n\geq 0$&  $\ZPdmjc=\left[{\scriptsize\begin{array}{cccc}D&\!\!\!0&\!\!\!B&\!\!\!C\\0&\!\!\!A&\!\!\!0&\!\!\!-C\end{array}}\right]$\sk\\
\hline\hline

$P(2n,1)$,& $\crk(\DML(\ZPdmj,\ZpPdmj,-D,n)$,\sk\\
$n\geq
0$&$\ZpPdmj=\left[{\scriptsize\begin{array}{ccc}B&\!\!\!C&\!\!\!D\end{array}}\right]$,
$\ZPdmj=\left[{\scriptsize\begin{array}{cccc}C&\!\!\!A&\!\!\!0&\!\!\!0\\-C&\!\!\!0&\!\!\!B&\!\!\!D\end{array}}\right]$\sk\\
\hline
$P(2n,2)$,& $\crk(\DML(\ZPdmd,\ZpPdmd,-A,n)$,\sk\\
$n\geq
0$&$\ZpPdmd=\left[{\scriptsize\begin{array}{ccc}C&\!\!\!D&\!\!\!A\end{array}}\right]$,
$\ZPdmd=\left[{\scriptsize\begin{array}{cccc}D&\!\!\!B&\!\!\!0&\!\!\!0\\-D&\!\!\!0&\!\!\!C&\!\!\!A\end{array}}\right]$\sk\\
\hline
$P(2n,3)$,& $\crk(\DML(\ZPdmt,\ZpPdmt,-B,n)$,\sk\\
$n\geq
0$&$\ZpPdmt=\left[{\scriptsize\begin{array}{ccc}D&\!\!\!A&\!\!\!B\end{array}}\right]$,
$\ZPdmt=\left[{\scriptsize\begin{array}{cccc}A&\!\!\!C&\!\!\!0&\!\!\!0\\-A&\!\!\!0&\!\!\!D&\!\!\!B\end{array}}\right]$\sk\\
\hline
$P(2n,4)$,& $\crk(\DML(\ZPdmc,\ZpPdmc,-C,n)$,\sk\\
$n\geq
0$&$\ZpPdmc=\left[{\scriptsize\begin{array}{ccc}A&\!\!\!B&\!\!\!C\end{array}}\right]$,
$\ZPdmc=\left[{\scriptsize\begin{array}{cccc}B&\!\!\!D&\!\!\!0&\!\!\!0\\-B&\!\!\!0&\!\!\!A&\!\!\!C\end{array}}\right]$\sk\\

\hline\hline $I(0,0)$&$n_0$\sk\\
\hline $I(n,0)$, &$\crk(\DM(\ZIm0,\ZpIm0,D\oplus C,n-1))$,\sk\\
$n\geq
1$&$\ZpIm0=\left[{\scriptsize\begin{array}{cccc}D&\!\!\!C&\!\!\!0&\!\!\!0\\0&\!\!\!-C&\!\!\!0&\!\!\!B\\-D&\!\!\!0&\!\!\!A&\!\!\!0\end{array}}\right]$,\
$\ZIm0=\left[{\scriptsize\begin{array}{cccc}0&\!\!\!-C&\!\!\!0&\!\!\!B\\-D&\!\!\!0&\!\!\!A&\!\!\!0\end{array}}\right]$\sk\\
\hline\hline

$I(2n+1,1)$, &$\crk(\DM(\ZIdmjj,A,-C,n))$,\sk\\
$n\geq 0$&
$\ZIdmjj=\left[{\scriptsize\begin{array}{cccc}C&\!\!\!D&\!\!\!0&\!\!\!B\\0&\!\!\!-D&\!\!\!A&\!\!\!0\end{array}}\right]$\sk\\

\hline $I(2n+1,2)$, &$\crk(\DM(\ZIdmjd,B,-D,n))$,\sk\\
$n\geq 0$&
$\ZIdmjd=\left[{\scriptsize\begin{array}{cccc}D&\!\!\!A&\!\!\!0&\!\!\!C\\0&\!\!\!-A&\!\!\!B&\!\!\!0\end{array}}\right]$\sk\\

\hline $I(2n+1,3)$, &$\crk(\DM(\ZIdmjt,C,-A,n))$,\sk\\
$n\geq 0$&
$\ZIdmjt=\left[{\scriptsize\begin{array}{cccc}A&\!\!\!B&\!\!\!0&\!\!\!D\\0&\!\!\!-B&\!\!\!C&\!\!\!0\end{array}}\right]$\sk\\

\hline $I(2n+1,4)$, &$\crk(\DM(\ZIdmjc,D,-B,n))$,\sk\\
$n\geq 0$&
$\ZIdmjc=\left[{\scriptsize\begin{array}{cccc}B&\!\!\!C&\!\!\!0&\!\!\!A\\0&\!\!\!-C&\!\!\!D&\!\!\!0\end{array}}\right]$\sk\\

\hline\hline
\end{tabular}

\begin{tabular}{ll}

\hline\hline $I(0,1)$&$n_1$\sk\\
\hline $I(2n,1)$, &$\crk(\DM(\ZIdmj,\ZpIdmj,-A,n-1))$,\sk\\
$n\geq
1$&$\ZpIdmj=\left[{\scriptsize\begin{array}{ccc}B&\!\!\!0&\!\!\!D\\0&\!\!\!C&\!\!\!-D\end{array}}\right]$,\
$\ZIdmj=\left[{\scriptsize\begin{array}{cccc}A&\!\!\!B&\!\!\!0&\!\!\!D\\0&\!\!\!0&\!\!\!C&\!\!\!-D\end{array}}\right]$\sk\\

\hline $I(0,2)$&$n_2$\sk\\
\hline $I(2n,2)$, &$\crk(\DM(\ZIdmd,\ZpIdmd,-B,n-1))$,\sk\\
$n\geq
1$&$\ZpIdmd=\left[{\scriptsize\begin{array}{ccc}C&\!\!\!0&\!\!\!A\\0&\!\!\!D&\!\!\!-A\end{array}}\right]$,\
$\ZIdmd=\left[{\scriptsize\begin{array}{cccc}B&\!\!\!C&\!\!\!0&\!\!\!A\\0&\!\!\!0&\!\!\!D&\!\!\!-A\end{array}}\right]$\sk\\

\hline $I(0,3)$&$n_3$\sk\\
\hline $I(2n,3)$, &$\crk(\DM(\ZIdmt,\ZpIdmt,-C,n-1))$,\sk\\
$n\geq
1$&$\ZpIdmt=\left[{\scriptsize\begin{array}{ccc}D&\!\!\!0&\!\!\!B\\0&\!\!\!A&\!\!\!-B\end{array}}\right]$,\
$\ZIdmt=\left[{\scriptsize\begin{array}{cccc}C&\!\!\!D&\!\!\!0&\!\!\!B\\0&\!\!\!0&\!\!\!A&\!\!\!-B\end{array}}\right]$\sk\\

\hline $I(0,4)$&$n_4$\sk\\
\hline $I(2n,4)$, &$\crk(\DM(\ZIdmc,\ZpIdmc,-D,n-1))$,\sk\\
$n\geq
1$&$\ZpIdmc=\left[{\scriptsize\begin{array}{ccc}A&\!\!\!0&\!\!\!C\\0&\!\!\!B&\!\!\!-C\end{array}}\right]$,\
$\ZIdmc=\left[{\scriptsize\begin{array}{cccc}D&\!\!\!A&\!\!\!0&\!\!\!C\\0&\!\!\!0&\!\!\!B&\!\!\!-C\end{array}}\right]$\sk\\
\hline\hline
 $R(l,\la)$, \ \ \ &
$\crk(\DM(\ZRlla,\ZRlla,-D,l-1)),$\sk\\
$l\geq 1,\la\neq
0,1$&$\ZRlla=\left[{\scriptsize\begin{array}{cccc}D&\!\!\!C&\!\!\!B&\!\!\!0\\-\la
D&\!\!\!-C&\!\!\!0&\!\!\!A\end{array}}\right]$\sk\\

\hline
\hline $R(0,2l,0)$,&$\crk(\DM(\ZRzdlz,\ZRzdlz,-D,l-1)),$\sk\\
$l\geq
1$&$\ZRzdlz=\left[{\scriptsize\begin{array}{cccc}D&\!\!\!C&\!\!\!B&\!\!\!0\\0&\!\!\!-C&\!\!\!0&\!\!\!A\end{array}}\right]$\sk\\
\hline $R(1,2l,0)$,&$\crk(\DM(\ZRjdlz,\ZRjdlz,-C,l-1)),$\sk\\
$l\geq
1$&$\ZRjdlz=\left[{\scriptsize\begin{array}{cccc}C&\!\!\!D&\!\!\!A&\!\!\!0\\0&\!\!\!-D&\!\!\!0&\!\!\!B\end{array}}\right]$\sk\\

\hline $R(0,2l,1)$,&$\crk(\DM(\ZRzdlj,\ZRzdlj,-B,l-1)),$\sk\\
$l\geq
1$&$\ZRzdlj=\left[{\scriptsize\begin{array}{cccc}B&\!\!\!D&\!\!\!A&\!\!\!0\\0&\!\!\!-D&\!\!\!0&\!\!\!C\end{array}}\right]$\sk\\
\hline $R(1,2l,1)$,&$\crk(\DM(\ZRjdlj,\ZRjdlj,-D,l-1)),$\sk\\
$l\geq
1$&$\ZRjdlj=\left[{\scriptsize\begin{array}{cccc}D&\!\!\!B&\!\!\!C&\!\!\!0\\0&\!\!\!-B&\!\!\!0&\!\!\!A\end{array}}\right]$\sk\\

\hline $R(0,2l,\infty)$,&$\crk(\DM(\ZRzdli,\ZRzdli,-D,l-1)),$\sk\\
$l\geq
1$&$\ZRzdli=\left[{\scriptsize\begin{array}{cccc}D&\!\!\!C&\!\!\!A&\!\!\!0\\0&\!\!\!-C&\!\!\!0&\!\!\!B\end{array}}\right]$\sk\\
\hline $R(1,2l,\infty)$,&$\crk(\DM(\ZRjdli,\ZRjdli,-C,l-1)),$\sk\\
$l\geq
1$&$\ZRjdli=\left[{\scriptsize\begin{array}{cccc}C&\!\!\!D&\!\!\!B&\!\!\!0\\0&\!\!\!-D&\!\!\!0&\!\!\!A\end{array}}\right]$\sk\\

\hline\hline

$R(0,2l-1,0)$, \ \ \ &
$\crk(\DM(\ZRzdljz,\ZpRzdljz,-B,l-1)),$\sk\\
$l\geq 1$&
$\ZpRzdljz=\left[{\scriptsize\begin{array}{cc}A&\!\!\!C\end{array}}\right]
$,\
$\ZRzdljz=\left[{\scriptsize\begin{array}{cccc}B&\!\!\!D&\!\!\!0&\!\!\!A\\0&\!\!\!0&\!\!\!C&\!\!\!-A\end{array}}\right]$\sk\\
\hline $R(1,2l-1,0)$, \ \ \ &
$\crk(\DM(\ZRjdljz,\ZpRjdljz,-A,l-1)),$\sk\\
$l\geq 1$&
$\ZpRjdljz=\left[{\scriptsize\begin{array}{cc}B&\!\!\!D\end{array}}\right]
$,\
$\ZRjdljz=\left[{\scriptsize\begin{array}{cccc}A&\!\!\!C&\!\!\!0&\!\!\!B\\0&\!\!\!0&\!\!\!D&\!\!\!-B\end{array}}\right]$\sk\\

\hline $R(0,2l-1,1)$, \ \ \ &
$\crk(\DM(\ZRzdljj,\ZpRzdljj,-A,l-1)),$\sk\\
$l\geq 1$&
$\ZpRzdljj=\left[{\scriptsize\begin{array}{cc}C&\!\!\!D\end{array}}\right]
$,\
$\ZRzdljj=\left[{\scriptsize\begin{array}{cccc}A&\!\!\!B&\!\!\!0&\!\!\!C\\0&\!\!\!0&\!\!\!D&\!\!\!-C\end{array}}\right]$\sk\\
\hline

\end{tabular}

\begin{tabular}{ll}

\hline $R(1,2l-1,1)$, \ \ \ &
$\crk(\DM(\ZRjdljj,\ZpRjdljj,-C,l-1)),$\sk\\
$l\geq 1$&
$\ZpRjdljj=\left[{\scriptsize\begin{array}{cc}A&\!\!\!B\end{array}}\right]
$,\
$\ZRjdljj=\left[{\scriptsize\begin{array}{cccc}C&\!\!\!D&\!\!\!0&\!\!\!A\\0&\!\!\!0&\!\!\!B&\!\!\!-A\end{array}}\right]$\sk\\

\hline $R(0,2l-1,\infty)$, \ \ \ &
$\crk(\DM(\ZRzdlji,\ZpRzdlji,-A,l-1)),$\sk\\
$l\geq 1$&
$\ZpRzdlji=\left[{\scriptsize\begin{array}{cc}B&\!\!\!C\end{array}}\right]
$,\
$\ZRzdlji=\left[{\scriptsize\begin{array}{cccc}A&\!\!\!D&\!\!\!0&\!\!\!B\\0&\!\!\!0&\!\!\!C&\!\!\!-B\end{array}}\right]$\sk\\
\hline $R(1,2l-1,\infty)$, \ \ \ &
$\crk(\DM(\ZRjdlji,\ZpRjdlji,-B,l-1)),$\sk\\
$l\geq 1$&
$\ZpRjdlji=\left[{\scriptsize\begin{array}{cc}A&\!\!\!D\end{array}}\right]
$,\
$\ZRjdlji=\left[{\scriptsize\begin{array}{cccc}B&\!\!\!C&\!\!\!0&\!\!\!A\\0&\!\!\!0&\!\!\!D&\!\!\!-A\end{array}}\right]$\sk\\

\hline\hline

\end{tabular}

\end{center}

 }

\bigskip

{\bf 2.} In the proof of Proposition 1 we essentially use the
matrix description of all the indecomposable modules over the four
subspace algebra $\La$, which we recall below (see \cite{Siks,SS};
note that in the description of indecomposables from the cited
handbooks one can find some number of misprints and little errors,
which sometimes caused that presented modules were decomposable;
we corrected the description here).

\smallskip For $m,n\geq 0$, by $0_{m\times n}\in\MM_{m\times n}$ or  $0_{n}\in\MM_{n\times n}$ we denote the zero
matrix and by $I_n\in\MM_{n\times n}$,  the identity matrix.
Moreover we set
\begin{center}
\begin{tabular}{rccclrcccl}
$\wlj_{n,n+1}$\!\!\!\!&$=$\!\!\!\!&$\left[{\scriptsize\begin{array}{ccccc}
\!\!\!\!1 & \!\!\!\! & \!\!\!\! &\!\!\!\!&\!\!\!\!
\\\!\!\!\! &
\!\!\!\! 1&\!\!\!\!& \!\!\!\!&\!\!\!\!
\\
\!\!\!\!&\!\!\!\!& \!\!\!\!\ddots & \!\!\!\!&
\!\!\!\!\\
\!\!\!\!\!\!\!\!& \!\!\!\!& \!\!\!\!&
\!\!\!\!1&\!\!\!\!0\!\!\!\!\end{array}}\right]$\!\!\!\!&$\in$\!\!\!\!&$\MM_{n\times
(n+1)}$,&$\wld_{n,n+1}$\!\!\!\!&$=$\!\!\!\!&$\left[{\scriptsize\begin{array}{ccccc}
\!\!\!\!0 & \!\!\!\! 1& \!\!\!\! &\!\!\!\!&\!\!\!\!
\\\!\!\!\! &
\!\!\!\! &\!\!\!\!\ddots& \!\!\!\!&\!\!\!\!
\\
\!\!\!\!&\!\!\!\!& \!\!\!\! & \!\!\!\!1&
\!\!\!\!\\
\!\!\!\!\!\!\!\!& \!\!\!\!& \!\!\!\!&
\!\!\!\!&\!\!\!\!1\!\!\!\!\end{array}}\right]$\!\!\!\!&$\in$\!\!\!\!&$\MM_{n\times
(n+1)}$,\smallskip
\\

$\Iod_n$\!\!\!\!&$=$\!\!\!\!&$\left[{\scriptsize\begin{array}{cccc}
\!\!\!\! & \!\!\!\! & \!\!\!\! &1\!\!\!\!\\\!\!\!\! & \!\!\!\!
&1\!\!\!\!& \!\!\!\!
\\ \!\!\!\!& \udots\!\!\!\! & \!\!\!\!& \!\!\!\!\\\!\!\!\!1\!\!\!\!& \!\!\!\!& \!\!\!\!& \!\!\!\!\end{array}}\right]$\!\!\!\!&$\in$\!\!\!\!&$\MM_{n\times
n}$,&$J_n(\la)$\!\!\!\!&$=$\!\!\!\!&$\left[{\scriptsize\begin{array}{cccccc}
\!\!\!\! \la& \!\!\!\! 1& \!\!\!\!
&\!\!\!\!&\!\!\!\!&\!\!\!\!\\
\!\!\!\! & \!\!\!\! \la&\!\!\!\!1& \!\!\!\!&\!\!\!\!&\!\!\!\!
\\
\!\!\!\!& \!\!\!\! & \!\!\!\!&
\!\!\!\!\ddots&\!\!\!\!&\!\!\!\!\\
\!\!\!\!\!\!\!\!& \!\!\!\!& \!\!\!\!&
\!\!\!\!&\!\!\!\!\la&\!\!\!\!1\\
\!\!\!\!&\!\!\!\!&\!\!\!\!&\!\!\!\!&\!\!\!\!&\!\!\!\!\la\!\!\!\!\end{array}}\right]$\!\!\!\!&$\in$\!\!\!\!&$\MM_{n\times
n}$,

\end{tabular}
\end{center}

\noindent for $n\geq 0$, $\la\in k$, where the remaining
coefficients are zero (note that $J_n(\la)$ is the
upper-triangular Jordan block with eigenvalue $\la$).

\medskip

{\bf Proposition.} {\em Let $\La$ be the four subspace algebra.
Then the fixed matrix representatives of all the isomorphism
classes of indecomposable $\La$-modules and their dimension
vectors are given in the following table: }

\begin{center}
\begin{tabular}{ll}
\hline\hline $P(n,0),$\ \ \ \ \ \ & $[\,2n+1,\,n,\,n,\,n,\,n\,],$\sk\\
$ n\geq
0$&$\left(\wmm{I_n}{0_{n+1,n}},\wmm{0_{n+1,n}}{I_n},\wmmm{0_{1,n}}{I_n}{\Iod_n},\wmmm{I_n}{\Iod_n}{0_{1,n}}\right)$\sk\\
\hline\hline $P(2n+1,1),$\ \ \ \ \ \ & $[\,2n+2,\,n,\,n+1,\,n+1,\,n+1\,],$\sk\\
$ n\geq
0$&$\left(\wmmmm{0_{1,n}}{I_n}{I_n}{0_{1,n}},\wmm{I_{n+1}}{0_{n+1}},\wmm{0_{n+1}}{I_{n+1}},\wmm{I_{n+1}}{I_{n+1}}\right)$\sk\\

\hline $P(2n+1,2),$\ \ \ \ \ \ & $[\,2n+2,\,n+1,\,n,\,n+1,\,n+1\,],$\sk\\
$ n\geq
0$&$\left(\wmm{I_{n+1}}{I_{n+1}},\wmmmm{0_{1,n}}{I_n}{I_n}{0_{1,n}},\wmm{I_{n+1}}{0_{n+1}},\wmm{0_{n+1}}{I_{n+1}}\right)$\sk\\

\hline $P(2n+1,3),$\ \ \ \ \ \ & $[\,2n+2,\,n+1,\,n+1,\,n,\,n+1\,],$\sk\\
$ n\geq
0$&$\left(\wmm{0_{n+1}}{I_{n+1}},\wmm{I_{n+1}}{I_{n+1}},\wmmmm{0_{1,n}}{I_n}{I_n}{0_{1,n}},\wmm{I_{n+1}}{0_{n+1}}\right)$\sk\\
\hline
\end{tabular}

\begin{tabular}{ll}

\hline $P(2n+1,4),$\ \ \ \ \ \ & $[\,2n+2,\,n+1,\,n+1,\,n+1,\,n]\,,$\sk\\
$ n\geq
0$&$\left(\wmm{I_{n+1}}{0_{n+1}},\wmm{0_{n+1}}{I_{n+1}},\wmm{I_{n+1}}{I_{n+1}},\wmmmm{0_{1,n}}{I_n}{I_n}{0_{1,n}}\right)$\sk\\

\hline \hline $P(2n,1),$\ \ \ \ \ \ & $[\,2n+1,\,n+1,\,n,\,n,\,n]\,,$\sk\\
$ n\geq
0$&$\left(\wmm{I_{n+1}}{0_{n,n+1}},\wmm{0_{n+1,n}}{I_{n}},\wmmm{0_{1,n}}{I_{n}}{I_{n}},\wmmm{I_n}{0_{1,n}}{I_n}\right)$\sk\\

\hline $P(2n,2),$\ \ \ \ \ \ & $[\,2n+1,\,n,\,n+1,\,n,\,n]$\,,\sk\\
$ n\geq
0$&$\left(\wmmm{I_n}{0_{1,n}}{I_n},\wmm{I_{n+1}}{0_{n,n+1}},\wmm{0_{n+1,n}}{I_{n}},\wmmm{0_{1,n}}{I_{n}}{I_{n}}\right)$\sk\\

\hline $P(2n,3),$\ \ \ \ \ \ & $[\,2n+1,\,n,\,n,\,n+1,\,n]$\,,\sk\\
$ n\geq
0$&$\left(\wmmm{0_{1,n}}{I_{n}}{I_{n}},\wmmm{I_n}{0_{1,n}}{I_n},\wmm{I_{n+1}}{0_{n,n+1}},\wmm{0_{n+1,n}}{I_{n}}\right)$\sk\\

\hline $P(2n,4),$\ \ \ \ \ \ & $[\,2n+1,\,n,\,n,\,n,\,n+1]$\,,\sk\\
$ n\geq
0$&$\left(\wmm{0_{n+1,n}}{I_{n}},\wmmm{0_{1,n}}{I_{n}}{I_{n}},\wmmm{I_n}{0_{1,n}}{I_n},\wmm{I_{n+1}}{0_{n,n+1}}\right)$\sk\\

\hline \hline $I(n,0),$\ \ \ \ \ \ & $[\,2n+1,\,n+1,\,n+1,\,n+1,\,n+1]\,,$\sk\\
$ n\geq
0$&$\left(\wmm{0_{n,n+1}}{I_{n+1}},\wmm{I_{n+1}}{0_{n,n+1}},\wmm{\Iod_{n+1}}{\wlj_{n,n+1}},\wmm{\wld_{n,n+1}}{\Iod_{n+1}}\right)$\sk\\

\hline \hline $I(2n+1,1),$\ \ \ \ \ \ & $[\,2n+1,\,n,\,n+1,\,n+1,\,n+1]\,,$\sk\\
$ n\geq
0$&$\left(\wmm{0_{n+1,n}}{I_{n}},\wmm{I_{n+1}}{0_{n,n+1}},\wmm{I_{n+1}}{\wld_{n,n+1}},\wmm{I_{n+1}}{\wlj_{n,n+1}}\right)$\sk\\

\hline $I(2n+1,2),$\ \ \ \ \ \ & $[\,2n+1,\,n+1,\,n,\,n+1,\,n+1]\,,$\sk\\
$ n\geq
0$&$\left(\wmm{I_{n+1}}{\wlj_{n,n+1}},\wmm{0_{n+1,n}}{I_{n}},\wmm{I_{n+1}}{0_{n,n+1}},\wmm{I_{n+1}}{\wld_{n,n+1}}\right)$\sk\\

\hline $I(2n+1,3),$\ \ \ \ \ \ & $[\,2n+1,\,n+1,\,n+1,\,n,\,n+1]\,,$\sk\\
$ n\geq
0$&$\left(\wmm{I_{n+1}}{\wld_{n,n+1}},\wmm{I_{n+1}}{\wlj_{n,n+1}},\wmm{0_{n+1,n}}{I_{n}},\wmm{I_{n+1}}{0_{n,n+1}}\right)$\sk\\

\hline $I(2n+1,4),$\ \ \ \ \ \ & $[\,2n+1,\,n+1,\,n+1,\,n+1,\,n]\,,$\sk\\
$ n\geq
0$&$\left(\wmm{I_{n+1}}{0_{n,n+1}},\wmm{I_{n+1}}{\wld_{n,n+1}},\wmm{I_{n+1}}{\wlj_{n,n+1}},\wmm{0_{n+1,n}}{I_{n}}\right)$\sk\\

\hline \hline $I(2n,1),$\ \ \ \ \ \ & $[\,2n,\,n+1,\,n,\,n,\,n]\,,$\sk\\
$ n\geq
0$&$\left(\wmm{\wlj_{n,n+1}}{\wld_{n,n+1}},\wmm{0_n}{I_n},\wmm{I_{n}}{0_n},\wmm{I_{n}}{I_n}\right)$\sk\\

\hline $I(2n,2),$\ \ \ \ \ \ & $[\,2n,\,n,\,n+1,\,n,\,n]\,,$\sk\\
$ n\geq
0$&$\left(\wmm{I_{n}}{I_n},\wmm{\wlj_{n,n+1}}{\wld_{n,n+1}},\wmm{0_n}{I_n},\wmm{I_{n}}{0_n}\right)$\sk\\

\hline $I(2n,3),$\ \ \ \ \ \ & $[\,2n,\,n,\,n,\,n+1,\,n]\,,$\sk\\
$ n\geq
0$&$\left(\wmm{I_{n}}{0_n},\wmm{I_{n}}{I_n},\wmm{\wlj_{n,n+1}}{\wld_{n,n+1}},\wmm{0_n}{I_n}\right)$\sk\\

\hline $I(2n,4),$\ \ \ \ \ \ & $[\,2n,\,n,\,n,\,n,\,n+1]\,,$\sk\\
$ n\geq
0$&$\left(\,\wmm{0_n}{I_n},\wmm{I_{n}}{0_n},\wmm{I_{n}}{I_n},\wmm{\wlj_{n,n+1}}{\wld_{n,n+1}}\,\right)$\sk\\
\hline

\end{tabular}

\begin{tabular}{ll}
\hline\hline $R(l,\la),$\ \ \ \ \ \ & $[\,2l,\,l,\,l,\,l,\,l]\,,$\sk\\
$ l\geq
1, \la\in k\setminus\{0,1\}$&$\left(\,\wmm{I_l}{0_l},\wmm{0_l}{I_l},\wmm{I_{l}}{I_l},\wmm{J_{l}(\la)}{I_l}\,\right)$\sk\\
\hline\hline
 $R(0,2l,0),$\ \ \ \ \ \ & $[\,2l,\,l,\,l,\,l,\,l]\,,$\sk\\
$ l\geq
1$&$\left(\,\wmm{I_l}{0_l},\wmm{0_l}{I_l},\wmm{I_{l}}{I_l},\wmm{J_{l}(0)}{I_l}\,\right)$\sk\\
\hline
 $R(1,2l,0),$\ \ \ \ \ \ & $[\,2l,\,l,\,l,\,l,\,l]\,,$\sk\\
$ l\geq
1$&$\left(\,\wmm{0_l}{I_l},\wmm{I_l}{0_l},\wmm{J_{l}(0)}{I_l},\wmm{I_l}{I_l}\,\right)$\sk\\
\hline
 $R(0,2l,1),$\ \ \ \ \ \ & $[\,2l,\,l,\,l,\,l,\,l]\,,$\sk\\
$ l\geq
1$&$\left(\,\wmm{0_l}{I_l},\wmm{J_l(0)}{I_l},\wmm{I_l}{0_l},\wmm{I_l}{I_l}\,\right)$\sk\\
\hline
 $R(1,2l,1),$\ \ \ \ \ \ & $[\,2l,\,l,\,l,\,l,\,l]\,,$\sk\\
$ l\geq
1$&$\left(\,\wmm{I_l}{0_l},\wmm{I_l}{I_l},\wmm{0_l}{I_l},\wmm{J_l(0)}{I_l}\,\right)$\sk\\
\hline
 $R(0,2l,\infty),$\ \ \ \ \ \ & $[\,2l,\,l,\,l,\,l,\,l]\,,$\sk\\
$ l\geq
1$&$\left(\,\wmm{0_l}{I_l},\wmm{I_l}{0_l},\wmm{I_l}{I_l},\wmm{J_l(0)}{I_l}\,\right)$\sk\\
\hline
 $R(1,2l,\infty),$\ \ \ \ \ \ & $[\,2l,\,l,\,l,\,l,\,l]\,,$\sk\\
$ l\geq
1$&$\left(\,\wmm{I_l}{0_l},\wmm{0_l}{I_l},\wmm{J_l(0)}{I_l},\wmm{I_l}{I_l}\,\right)$\sk\\
\hline\hline
$R(0,2l-1,0),$\ \ \ \ \ \ & $[\,2l-1,\,l-1,\,l,\,l-1,\,l]\,,$\sk\\
$ l\geq
1$&$\left(\,\wmmm{I_{l-1}}{0_{1,l-1}}{I_{l-1}},\wmm{\wlj_{l-1,l}}{I_l},\wmm{I_{l-1}}{0_{l,l-1}},\wmm{0_{l-1,l}}{I_l}\,\right)$\sk\\
\hline

$R(1,2l-1,0),$\ \ \ \ \ \ & $[\,2l-1,\,l,\,l-1,\,l,\,l-1]\,,$\sk\\
$ l\geq
1$&$\left(\,\wmm{\wlj_{l-1,l}}{I_l},\wmmm{I_{l-1}}{0_{1,l-1}}{I_{l-1}},\wmm{0_{l-1,l}}{I_l},\wmm{I_{l-1}}{0_{l,l-1}}\,\right)$\sk\\
\hline

$R(0,2l-1,1),$\ \ \ \ \ \ & $[\,2l-1,\,l,\,l,\,l-1,\,l-1]\,,$\sk\\
$ l\geq
1$&$\left(\,\wmm{\wlj_{l-1,l}}{I_l},\wmm{0_{l-1,l}}{I_l},\wmmm{I_{l-1}}{0_{1,l-1}}{I_{l-1}},\wmm{I_{l-1}}{0_{l,l-1}}\,\right)$\sk\\
\hline
$R(1,2l-1,1),$\ \ \ \ \ \ & $[\,2l-1,\,l-1,\,l-1,\,l,\,l]\,,$\sk\\
$ l\geq
1$&$\left(\,\wmmm{I_{l-1}}{0_{1,l-1}}{I_{l-1}},\wmm{I_{l-1}}{0_{l,l-1}},\wmm{\wlj_{l-1,l}}{I_l},\wmm{0_{l-1,l}}{I_l}\,\right)$\sk\\
\hline

$R(0,2l-1,\infty),$\ \ \ \ \ \ & $[\,2l-1,\,l,\,l-1,\,l-1,\,l]\,,$\sk\\
$ l\geq
1$&$\left(\,\wmm{\wlj_{l-1,l}}{I_l},\wmmm{I_{l-1}}{0_{1,l-1}}{I_{l-1}},\wmm{I_{l-1}}{0_{l,l-1}},\wmm{0_{l-1,l}}{I_l}\,\right)$\sk\\
\hline

$R(1,2l-1,\infty),$\ \ \ \ \ \ & $[\,2l-1,\,l-1,\,l,\,l,\,l-1]\,,$\sk\\
$ l\geq
1$&$\left(\,\wmmm{I_{l-1}}{0_{1,l-1}}{I_{l-1}},\wmm{\wlj_{l-1,l}}{I_l},\wmm{0_{l-1,l}}{I_l},\wmm{I_{l-1}}{0_{l,l-1}}\,\right)$\sk\\
\hline \hline

\end{tabular}

\end{center}

\bigskip

One proves the proposition using the
standard arguments as in \cite{SS, Siks}.




\bigskip

{\bf Remark.} Note that for any $l\geq 1$ we have
$R(0,2l,0)=R(l,\la)$ if we substitute $\la:=0$. One can also
easily check that $R(1,2l,1)$ is isomorphic to $R(l,\la)$ with
$\la:=1$. Therefore for every $\la\in k$  (not only for $\la\in
k\setminus\{0,1\}$), the tube $\CT_\la$ contains the
indecomposable module
 $R(l,\la)$, for any $l\geq 1$.

\bigskip

 {\bf 3.} {\em Proof of Proposition 1.} Let $M=(A,B,C,D)$ be a fixed $\La$-module,
with $A\in\MM_{n_0\times n_1}$, $B\in\MM_{n_0\times n_2}$,
$C\in\MM_{n_0\times n_3}$, $D\in\MM_{n_0\times n_4}$. Recall
\cite{ASS} that for $\La$-module $X=(A',B',C',D')$ with
$A'\in\MM_{m_0\times m_1}$, $B'\in\MM_{m_0\times m_2}$,
$C'\in\MM_{m_0\times m_3}$, $D'\in\MM_{m_0\times m_4}$, any
$\La$-homomorphism $F:M\to X$ is  the collection
$(F_0,F_1,F_2,F_3,F_4)$ of $k$-homomorphisms $F_i:k^{n_i}\to
k^{m_i}$, $i\in Q_0$, satisfying the commutativity relations
$$F_0A=A'F_1,\ \ F_0B=B'F_2,\ \ F_0C=C'F_3,\ \
F_0D=D'F_4\leqno{(\ast)}$$ (we identify $F_i$ with its matrix from
$\MM_{m_i\times n_i}$ given in the standard bases). So the
dimension $[M,X]$ equals the dimension of the solution space of
the system $(\ast)$ of linear equations with the matrices $F_i$
treated as variables. We calculate below these dimensions for
every indecomposable $\La$-module $X$. We use the description of
indecomposables from Proposition 2.


\smallskip

\underline{Case 1: $X=P(n,0)$}. $P(0,0)$ is the projective
$\La$-module $P(0)=(\trivm_{1,0},\trivm_{1,0},\trivm_{1,0},$ $
\trivm_{1,0})$. So every homomorphism from $M$ to $P(0,0)$ is
given by the collection $(F_0,F_1,F_2,$ $F_3,F_4)=(y,0,0,0,0)$,
$y\in\MM_{1\times n_0}$, such that $yA=0$, $yB=0$, $yC=0$, $yD=0$
that is, satisfying the matrix equation $y[A\ B\ C\ D]=0$. Thus
$[M,P(0,0)]=\crk([A\ B\ C\ D])$.

Now fix $n\geq 1$. Recall that
$P(n,0)=\left(\wmm{I_n}{0_{n+1,n}},\wmm{0_{n+1,n}}{I_n},\wmmm{0_{1,n}}{I_n}{\Iod_n},\wmmm{I_n}{\Iod_n}{0_{1,n}}\right)$
(see Proposition 2). So every homomorphism from $M$ to $P(n,0)$ is
given by the collection of matrices $(F_0,F_1,F_2,F_3,F_4)\in
\MM_{(2n+1)\times n_0}\times\MM_{n\times n_1}\times\MM_{n\times
n_2}\times\MM_{n\times n_3}\times\MM_{n\times n_4}$ satisfying the
commutativity relations $$F_0A=\wmm{I_n}{0_{n+1,n}}F_1,\ \
F_0B=\wmm{0_{n+1,n}}{I_n}F_2,\ \
F_0C=\wmmm{0_{1,n}}{I_n}{\Iod_n}F_3,\ \
F_0D=\wmmm{I_n}{\Iod_n}{0_{1,n}}F_4.\leqno{(\ast)^{1}}$$
 Let us denote by $y_1,\ldots,y_{2n+1}$, $s_1,\ldots, s_n$,
$t_1,\ldots,t_n$, $u_1,\ldots,u_n$ and respectively
$w_1,\ldots,w_n$, the consecutive rows of matrices
$F_0,F_1,F_2,F_3,F_4$. So the system $(\ast)^{1}$ can be presented
in the form
$${\scriptsize \left\{\begin{array}{lll}
y_1A&\!\!\!\!=&\!\!\!\!s_1\\
&\!\!\!\!\vdots&\\
y_nA&\!\!\!\!=&\!\!\!\!s_n\\
y_{n+1}A&\!\!\!\!=&\!\!\!\!0\\
&\!\!\!\!\vdots&\\
y_{2n}A&\!\!\!\!=&\!\!\!\!0\\
y_{2n+1}A&\!\!\!\!=&\!\!\!\!0\\
\end{array}\right.\  \ \
\left\{\begin{array}{lll}
y_1B&\!\!\!\!=&\!\!\!\!0\\
&\!\!\!\!\vdots&\\
y_nB&\!\!\!\!=&\!\!\!\!0\\
y_{n+1}B&\!\!\!\!=&\!\!\!\!0\\
y_{n+2}B&\!\!\!\!=&\!\!\!\!t_1\\
&\!\!\!\!\vdots&\\
y_{2n+1}B&\!\!\!\!=&\!\!\!\!t_n
\end{array}\right.\ \ \
\left\{\begin{array}{lll}
y_1C&\!\!\!\!=&\!\!\!\!0\\
y_2C&\!\!\!\!=&\!\!\!\!u_1\\
&\!\!\!\!\vdots&\\
y_{n+1}C&\!\!\!\!=&\!\!\!\!u_n\\
y_{n+2}C&\!\!\!\!=&\!\!\!\!u_n\\
&\!\!\!\!\vdots&\\
y_{2n+1}C&\!\!\!\!=&\!\!\!\!u_1
\end{array}\right.\ \ \
\left\{\begin{array}{lll}
y_1D&\!\!\!\!=&\!\!\!\!w_1\\
&\!\!\!\!\vdots&\\
y_{n}D&\!\!\!\!=&\!\!\!\!w_n\\
y_{n+1}D&\!\!\!\!=&\!\!\!\!w_n\\
&\!\!\!\!\vdots&\\
y_{2n}D&\!\!\!\!=&\!\!\!\!w_1\\
y_{2n+1}D&\!\!\!\!=&\!\!\!\!0\\
\end{array}\right.
}\leqno{(\ast\ast)^1}
$$
Observe that  all the variables $s_i, t_i, u_i, w_i$ are
determined by $y_1,\ldots,$ $y_{2n+1}$, so clearly the dimensions
of the solution spaces for $(\ast\ast)^1$ and
$${\scriptsize \left\{\begin{array}{lll}
y_{n+1}A&\!\!\!\!=&\!\!\!\!0\\
y_{n+2}A&\!\!\!\!=&\!\!\!\!0\\
&\!\!\!\!\vdots&\\
y_{2n+1}A&\!\!\!\!=&\!\!\!\!0\\
\end{array}\right.\  \ \
\left\{\begin{array}{lll}
y_1B&\!\!\!\!=&\!\!\!\!0\\
y_2B&\!\!\!\!=&\!\!\!\!0\\
&\!\!\!\!\vdots&\\
y_{n+1}B&\!\!\!\!=&\!\!\!\!0\\
\end{array}\right.\ \ \
\left\{\begin{array}{lll}
y_1C&\!\!\!\!=&\!\!\!\!0\\
y_2C-y_{2n+1}C&\!\!\!\!=&\!\!\!\!0\\
y_3C-y_{2n}C&\!\!\!\!=&\!\!\!\!0\\
&\!\!\!\!\vdots&\\
y_{n+1}C-y_{n+2}C&\!\!\!\!=&\!\!\!\!0
\end{array}\right.\ \ \
\left\{\begin{array}{lll}
y_1D-y_{2n}D&\!\!\!\!=&\!\!\!\!0\\
y_{2}D-y_{2n-1}D&\!\!\!\!=&\!\!\!\!0\\
&\!\!\!\!\vdots&\\
y_{n}D-y_{n+1}D&\!\!\!\!=&\!\!\!\!0\\
y_{2n+1}D&\!\!\!\!=&\!\!\!\!0\\
\end{array}\right.
}\leqno{(\ast\!\ast\!\ast)^1}
$$

\noindent are equal.

For the benefit of the reader we discuss in details  the cases
$n=1$ and $2$. For $n=1$, the system $(\ast\!\ast\!\ast)^1$ has
the form
$${\scriptsize \left\{\begin{array}{lll}
y_{2}A&\!\!\!\!=&\!\!\!\!0\\
y_{3}A&\!\!\!\!=&\!\!\!\!0\\
\end{array}\right.\  \ \
\left\{\begin{array}{lll}
y_1B&\!\!\!\!=&\!\!\!\!0\\
y_2B&\!\!\!\!=&\!\!\!\!0\\
\end{array}\right.\ \ \
\left\{\begin{array}{lll}
y_1C&\!\!\!\!=&\!\!\!\!0\\
y_2C-y_{3}C&\!\!\!\!=&\!\!\!\!0\\
\end{array}\right.\ \ \
\left\{\begin{array}{lll}
y_1D-y_{2}D&\!\!\!\!=&\!\!\!\!0\\
y_{3}D&\!\!\!\!=&\!\!\!\!0\\
\end{array}\right.
}
$$
and we can equivalently write it in the following
matrix form
$$[y_2\ y_1\ y_3]\left[{\scriptsize\begin{array}{cccccccc}A&\!\!\!0&\!\!\!B&\!\!\!0&\!\!\!C&\!\!\!D&\!\!\!0&\!\!\!0\\0&\!\!\!0&\!\!\!0&\!\!\!B&\!\!\!0&\!\!\!-D&\!\!\!C&\!\!\!0\\
0&\!\!\!A&\!\!\!0&\!\!\!0&\!\!\!-C&\!\!\!0&\!\!\!0&\!\!\!D\\\end{array}}\right]=0$$
Note that the matrix appearing in the above equation equals
$$\CN_1=\DN(Z_{P(1,0)},Z'_{P(1,0)},C\oplus D,0),$$ so the dimension of the
solution space of this equation, which is equal to $[M,P(1,0)]$,
is $\cor(\CN_1)$, and this is precisely the assertion for $n=1$.

For $n=2$, the system $(\ast\!\ast\!\ast)^1$ has the form
$${\scriptsize \left\{\begin{array}{lll}
y_{3}A&\!\!\!\!=&\!\!\!\!0\\
y_{4}A&\!\!\!\!=&\!\!\!\!0\\
y_{5}A&\!\!\!\!=&\!\!\!\!0\\
\end{array}\right.\  \ \
\left\{\begin{array}{lll}
y_1B&\!\!\!\!=&\!\!\!\!0\\
y_2B&\!\!\!\!=&\!\!\!\!0\\
y_3B&\!\!\!\!=&\!\!\!\!0\\
\end{array}\right.\ \ \
\left\{\begin{array}{lll}
y_1C&\!\!\!\!=&\!\!\!\!0\\
y_2C-y_{5}C&\!\!\!\!=&\!\!\!\!0\\
y_3C-y_{4}C&\!\!\!\!=&\!\!\!\!0\\
\end{array}\right.\ \ \
\left\{\begin{array}{lll}
y_1D-y_{4}D&\!\!\!\!=&\!\!\!\!0\\
y_2D-y_{3}D&\!\!\!\!=&\!\!\!\!0\\
y_{5}D&\!\!\!\!=&\!\!\!\!0\\
\end{array}\right.
}
$$
and we can equivalently write it in the following matrix form
$$[y_3\ y_2\ y_4\ y_1\ y_5]\left[{\scriptsize\begin{array}{cccccc|cccc|cc}A&\!\!\!0&\!\!\!B&\!\!\!0&\!\!\!C&\!\!\!D&\!\!\!0&\!\!\!0&\!\!\!0&\!\!\!0&\!\!\!0&\!\!\!0\\0&\!\!\!0&\!\!\!0&\!\!\!B&\!\!\!0&\!\!\!-D&\!\!\!C&\!\!\!0&\!\!\!0&\!\!\!0&\!\!\!0&\!\!\!0\\
0&\!\!\!A&\!\!\!0&\!\!\!0&\!\!\!-C&\!\!\!0&\!\!\!0&\!\!\!D&\!\!\!0&\!\!\!0&\!\!\!0&\!\!\!0\\\hline
0&\!\!\!0&\!\!\!0&\!\!\!0&\!\!\!0&\!\!\!0&\!\!\!0&\!\!\!-D&\!\!\!0&\!\!\!B&\!\!\!C&\!\!\!0\\
0&\!\!\!0&\!\!\!0&\!\!\!0&\!\!\!0&\!\!\!0&\!\!\!-C&\!\!\!0&\!\!\!A&\!\!\!0&\!\!\!0&\!\!\!D\end{array}}\right]=0$$
Note that the matrix appearing in the above equation equals
$$\CN_2=\DN(Z_{P(2,0)},Z'_{P(2,0)},C\oplus D,1),$$ so the dimension of the
solution space of this equation, which is equal to $[M,P(2,0)]$,
is $\cor(\CN_2)$, and this is precisely the assertion for $n=2$.

Similarly one checks that for the general case $n\geq 1$, the
system $(\ast\!\ast\!\ast)^1$ is equivalent to the matrix equation
$$
[y_{n+1}\ y_n\ y_{n+2}\ y_{n-1}\ y_{n+3}\ y_{n-2}\ \ldots\ y_{2n}\
y_1\ y_{2n+1}]\cdot\CN_n=0,
$$
where $\CN_n=\DN(\ZPm0,\ZpPm0,C\oplus D,n-1)$. Therefore the
dimension $[M,P(n,0)]$ equals $\crk(\CN_n)$ and we are done.

\smallskip

\underline{Case 2: $X=P(2n+1,1)$}. Fix $n\geq 0$. Then every
homomorphism from $M$ to $P(2n+1,1)$ is given by the collection of
matrices $(F_0,F_1,F_2,F_3,F_4)\in \MM_{(2n+2)\times
n_0}\times\MM_{n\times n_1}\times\MM_{(n+1)\times
n_2}\times\MM_{(n+1)\times n_3}\times\MM_{(n+1)\times n_4}$
satisfying the commutativity relations
$$F_0A=\wmmmm{0_{1,n}}{I_n}{I_n}{0_{1,n}}F_1,\ \ F_0B=\wmm{I_{n+1}}{0_{n+1}}F_2,\
\ F_0C=\wmm{0_{n+1}}{I_{n+1}}F_3,\ \
F_0D=\wmm{I_{n+1}}{I_{n+1}}F_4\leqno{(\ast)^{2}}$$ (see $(\ast)$
and Proposition 2).
 Let us denote by $y_1,\ldots,y_{2n+2}$, $s_1,\ldots, s_n$,
$t_1,\ldots,t_{n+1}$, $u_1,\ldots,u_{n+1}$ and respectively
$w_1,\ldots,w_{n+1}$, the consecutive rows of matrices $F_0,F_1,$
$ F_2,F_3,F_4$. Similarly as in the case 1, the variables $s_i,
t_i, u_i, w_i$ are determined by $y_1,\ldots,$ $y_{2n+2}$, so
clearly the dimensions of the solution spaces for $(\ast)^2$ and

$${\scriptsize \left\{\begin{array}{lll}
y_{1}A&\!\!\!\!=&\!\!\!\!0\\
y_{2}A-y_{n+2}A&\!\!\!\!=&\!\!\!\!0\\
y_{3}A-y_{n+3}A&\!\!\!\!=&\!\!\!\!0\\
&\!\!\!\!\vdots&\\
y_{n+1}A-y_{2n+1}A&\!\!\!\!=&\!\!\!\!0\\
y_{2n+2}A&\!\!\!\!=&\!\!\!\!0\\
\end{array}\right.\  \ \
\left\{\begin{array}{lll}
y_{n+2}B&\!\!\!\!=&\!\!\!\!0\\
y_{n+3}B&\!\!\!\!=&\!\!\!\!0\\
&\!\!\!\!\vdots&\\
y_{2n+2}B&\!\!\!\!=&\!\!\!\!0\\
\end{array}\right.\ \ \
\left\{\begin{array}{lll}
y_1C&\!\!\!\!=&\!\!\!\!0\\
y_2C&\!\!\!\!=&\!\!\!\!0\\
&\!\!\!\!\vdots&\\
y_{n+1}C&\!\!\!\!=&\!\!\!\!0
\end{array}\right.\ \ \
\left\{\begin{array}{lll}
y_1D-y_{n+2}D&\!\!\!\!=&\!\!\!\!0\\
y_{2}D-y_{n+3}D&\!\!\!\!=&\!\!\!\!0\\
&\!\!\!\!\vdots&\\
y_{n+1}D-y_{2n+2}D&\!\!\!\!=&\!\!\!\!0\\
\end{array}\right.
}\leqno{(\ast\!\ast)^2}
$$

\noindent are equal. One checks that  the system $(\ast\!\ast)^2$
is equivalent to the matrix equation
$$
[y_{1}\ y_{n+2}\ y_{2}\ y_{n+3}\ y_{3}\ y_{n+4}\ \ldots\ y_{n}\
y_{2n+1}\ y_{n+1}\ y_{2n+2}]\cdot\CN_n=0,
$$
where $\CN_n=\DN(\ZP2m11,\ZP2m11,-A,n)$. Therefore the dimension
$[M,P(2n+1,1)]$ equals $\crk(\CN_n)$.

\smallskip
\underline{Case 3: $X=P(2n+1,i)$, $i\in\{2,3,4\}$}. The assertions
for this case follow immediately from the case 2, since the module
$P(2n+1,i)$, for $i\in\{2,3,4\}$, can be obtained from
$P(2n+1,i-1)$ after the cyclic permutation $(1,2,3,4)$ of vertices
of the quiver $Q$ (see Proposition 2).

\smallskip

\underline{Case 4: $X=P(2n,1)$}. Fix $n\geq 0$. Then every
homomorphism from $M$ to $P(2n,1)$ is given by the collection of
matrices $(F_0,F_1,F_2,F_3,F_4)\in \MM_{(2n+1)\times
n_0}\times\MM_{(n+1)\times n_1}\times\MM_{n\times
n_2}\times\MM_{n\times n_3}\times\MM_{n\times n_4}$ satisfying the
commutativity relations
$$F_0A=\wmm{I_{n+1}}{0_{n,n+1}}F_1,\ \ F_0B=\wmm{0_{n+1,n}}{I_{n}}F_2,\
\ F_0C=\wmmm{0_{1,n}}{I_{n}}{I_{n}}F_3,\ \
F_0D=\wmmm{I_n}{0_{1,n}}{I_n}F_4\leqno{(\ast)^{4}}$$ (see $(\ast)$
and Proposition 2).
 Let us denote by $y_1,\ldots,y_{2n+1}$, $s_1,\ldots, s_{n+1}$,
$t_1,\ldots,t_{n}$, $u_1,\ldots,u_{n}$ and respectively
$w_1,\ldots,w_{n}$, the consecutive rows of matrices $F_0,F_1,$ $
F_2,F_3,F_4$. Similarly as before, the variables $s_i, t_i, u_i,
w_i$ are determined by $y_1,\ldots,$ $y_{2n+1}$, so clearly the
dimensions of the solution spaces for $(\ast)^4$ and
$${\scriptsize \left\{\begin{array}{lll}
y_{n+2}A&\!\!\!\!=&\!\!\!\!0\\
y_{n+3}A&\!\!\!\!=&\!\!\!\!0\\
&\!\!\!\!\vdots&\\
y_{2n+1}A&\!\!\!\!=&\!\!\!\!0\\
\end{array}\right.\  \ \
\left\{\begin{array}{lll}
y_{1}B&\!\!\!\!=&\!\!\!\!0\\
y_{2}B&\!\!\!\!=&\!\!\!\!0\\
&\!\!\!\!\vdots&\\
y_{n+1}B&\!\!\!\!=&\!\!\!\!0\\
\end{array}\right.\ \ \
\left\{\begin{array}{lll}
y_1C&\!\!\!\!=&\!\!\!\!0\\
y_2C-y_{n+2}C&\!\!\!\!=&\!\!\!\!0\\
y_3C-y_{n+3}C&\!\!\!\!=&\!\!\!\!0\\
&\!\!\!\!\vdots&\\
y_{n+1}C-y_{2n+1}C&\!\!\!\!=&\!\!\!\!0
\end{array}\right.\ \ \
\left\{\begin{array}{lll}
y_1D-y_{n+2}D&\!\!\!\!=&\!\!\!\!0\\
y_{2}D-y_{n+3}D&\!\!\!\!=&\!\!\!\!0\\
&\!\!\!\!\vdots&\\
y_{n}D-y_{2n+1}D&\!\!\!\!=&\!\!\!\!0\\
y_{n+1}D&\!\!\!\!=&\!\!\!\!0\\
\end{array}\right.
}\leqno{(\ast\!\ast)^4}
$$

\noindent are equal.

For $n=0$, the system $(\ast\!\ast)^4$ has the form $y_1B=0$,
$y_1C=0$ and $y_1D=0$, and we can write it in the matrix form as
$y_1\CN_0=0$, where $\CN_0=[B\,C\,D]$. So the dimension
$[M,P(0,1)]$ equals $\cor(\CN_0)$ and this is the assertion for
$n=0$, since $\CN_0=\DML(Z_{P(0,1)},Z'_{P(0,1)},-D,0)$.

Similarly one checks that in the general case $n\geq 0$, the
system $(\ast\!\ast)^4$ is equivalent to the matrix equation
$$
[y_{1}\ y_{n+2}\ y_{2}\ y_{n+3}\ y_{3}\ y_{n+4}\ \ldots\ y_{n}\
y_{2n+1}\ y_{n+1}]\cdot\CN_n=0,
$$
where $\CN_n=\DML(\ZPdmj,\ZpPdmj,-D,n)$. Therefore the dimension
$[M,P(2n,1)]$ equals $\crk(\CN_n)$.

\smallskip

\underline{Case 5: $X=P(2n,i)$, $i\in\{2,3,4\}$}. The assertions
for this case follow immediately from the case 4, since the module
$P(2n,i)$, for $i\in\{2,3,4\}$, can be obtained from $P(2n,i-1)$
after the cyclic permutation $(1,2,3,4)$ of vertices of the quiver
$Q$ (see Proposition 2).

\smallskip

\underline{Case 6: $X=I(n,0)$}. $I(0,0)$ is the injective
$\La$-module
$I(0)=(\,\scriptsize{[1]},\scriptsize{[1]},\scriptsize{[1]},\scriptsize{[1]}\,)$.
So every homomorphism from $M$ to $I(0,0)$ is given by the
collection $(y,s,t,u,w)\in\MM_{1\times n_0}\times \MM_{1\times
n_1}\times\MM_{1\times n_2}\times\MM_{1\times
n_3}\times\MM_{1\times n_4}$, such that $yA=s$, $yB=t$, $yC=u$,
$yD=w$ that is, the variables $s,t,u,w$ are determined by $y$ and
$y$ can have an arbitrary value. Thus $[M,I(0,0)]={\rm
dim}(\MM_{1\times n_0})=n_0$.

For $n\geq 1$, every homomorphism $F:M\to I(n,0)$ is given by the
collection $(F_0,F_1,F_2,F_3,F_4)\in \MM_{(2n+1)\times
n_0}\times\MM_{(n+1)\times n_1}\times\MM_{(n+1)\times
n_2}\times\MM_{(n+1)\times n_3}\times\MM_{(n+1)\times n_4}$
satisfying the commutativity relations
$$F_0A=\wmm{0_{n,n+1}}{I_{n+1}}F_1,\ \ F_0B=\wmm{I_{n+1}}{0_{n,n+1}}F_2,\
\ F_0C=\wmm{\Iod_{n+1}}{\wlj_{n,n+1}}F_3,\ \
F_0D=\wmm{\wld_{n,n+1}}{\Iod_{n+1}}F_4\leqno{(\ast)^{6}}$$ (see
$(\ast)$ and Proposition 2).
 Let us denote by $y_1,\ldots,y_{2n+1}$, $s_1,\ldots, s_{n+1}$,
$t_1,\ldots,t_{n+1}$, $u_1,\ldots,u_{n+1}$ and respectively
$w_1,\ldots,w_{n+1}$, the consecutive rows of matrices $F_0,F_1,$
$ F_2,F_3,F_4$. Similarly as before, the variables $s_i, t_i, u_i,
w_i$ are determined by $y_1,\ldots,$ $y_{2n+1}$, so clearly the
dimensions of the solution spaces for $(\ast)^6$ and
$${\scriptsize \left\{\begin{array}{lll}
y_{1}A&\!\!\!\!=&\!\!\!\!0\\
y_{2}A&\!\!\!\!=&\!\!\!\!0\\
&\!\!\!\!\vdots&\\
y_{n}A&\!\!\!\!=&\!\!\!\!0\\
\end{array}\right.\  \ \
\left\{\begin{array}{lll}
y_{n+2}B&\!\!\!\!=&\!\!\!\!0\\
y_{n+3}B&\!\!\!\!=&\!\!\!\!0\\
&\!\!\!\!\vdots&\\
y_{2n+1}B&\!\!\!\!=&\!\!\!\!0\\
\end{array}\right.\ \ \
\left\{\begin{array}{lll}
y_2C-y_{2n+1}C&\!\!\!\!=&\!\!\!\!0\\
y_3C-y_{2n}C&\!\!\!\!=&\!\!\!\!0\\
&\!\!\!\!\vdots&\\
y_{n}C-y_{n+3}C&\!\!\!\!=&\!\!\!\!0\\
y_{n+1}C-y_{n+2}C&\!\!\!\!=&\!\!\!\!0
\end{array}\right.\ \ \
\left\{\begin{array}{lll}
y_1D-y_{2n}D&\!\!\!\!=&\!\!\!\!0\\
y_{2}D-y_{2n-1}D&\!\!\!\!=&\!\!\!\!0\\
&\!\!\!\!\vdots&\\
y_{n-1}D-y_{n+2}D&\!\!\!\!=&\!\!\!\!0\\
y_{n}D-y_{n+1}D&\!\!\!\!=&\!\!\!\!0\\
\end{array}\right.
}\leqno{(\ast\!\ast)^6}
$$

\noindent are equal.

One checks that the system $(\ast\!\ast)^6$ is equivalent to the
matrix equation
$$
[y_{n+1}\ y_{n+2}\ y_n\ y_{n+3}\ y_{n-1}\ y_{n+4}\ \ldots\ y_2\
y_{2n+1}\ y_1]\cdot\CN_n=0,
$$
where $\CN_n=\DM(\ZIm0,\ZpIm0,D\oplus C,n-1)$. Therefore the
dimension $[M,I(n,0)]$ equals $\crk(\CN_n)$.

\smallskip

\underline{Case 7: $X=I(2n+1, 1)$}. $I(1,1)$ has the form
$(\trivm_{1,0},\scriptsize{[1]},\scriptsize{[1]},\scriptsize{[1]}\,)$.
So every homomorphism $F:M\to I(1,1)$ is given by the collection
$(y,0,t,u,w)\in\MM_{1\times n_0}\times \MM_{0\times
n_1}\times\MM_{1\times n_2}\times\MM_{1\times
n_3}\times\MM_{1\times n_4}$, such that $yA=0$, $yB=t$, $yC=u$,
$yD=w$. Therefore the variables $t,u,w$ are determined by $y$, so
the dimension $[M,I(1,1)]$ equals $\crk(A)$ and this is precisely
the assertion for $n=0$ since $A=\DM(Z_{I(1,1)},A,-C,0)$.

For $n\geq 1$, every homomorphism $F:M\to I(2n+1,1)$ is given by
the collection $(F_0,F_1,F_2,F_3,F_4)\in \MM_{(2n+1)\times
n_0}\times\MM_{n\times n_1}\times\MM_{(n+1)\times
n_2}\times\MM_{(n+1)\times n_3}\times\MM_{(n+1)\times n_4}$
satisfying the commutativity relations
$$F_0A=\wmm{0_{n+1,n}}{I_{n}}F_1,\ \ F_0B=\wmm{I_{n+1}}{0_{n,n+1}}F_2,\
\ F_0C=\wmm{I_{n+1}}{\wld_{n,n+1}}F_3,\ \
F_0D=\wmm{I_{n+1}}{\wlj_{n,n+1}}F_4\leqno{(\ast)^{7}}$$ (see
$(\ast)$ and Proposition 2).
 Let us denote by $y_1,\ldots,y_{2n+1}$, $s_1,\ldots, s_{n}$,
$t_1,\ldots,t_{n+1}$, $u_1,\ldots,u_{n+1}$ and respectively
$w_1,\ldots,w_{n+1}$, the consecutive rows of matrices $F_0,F_1,$
$ F_2,F_3,F_4$. Similarly as before, the variables $s_i, t_i, u_i,
w_i$ are determined by $y_1,\ldots,$ $y_{2n+1}$, so clearly the
dimensions of the solution spaces for $(\ast)^7$ and
$${\scriptsize \left\{\begin{array}{lll}
y_{1}A&\!\!\!\!=&\!\!\!\!0\\
&\!\!\!\!\vdots&\\
y_{n}A&\!\!\!\!=&\!\!\!\!0\\
y_{n+1}A&\!\!\!\!=&\!\!\!\!0
\end{array}\right.\  \ \
\left\{\begin{array}{lll}
y_{n+2}B&\!\!\!\!=&\!\!\!\!0\\
y_{n+3}B&\!\!\!\!=&\!\!\!\!0\\
&\!\!\!\!\vdots&\\
y_{2n+1}B&\!\!\!\!=&\!\!\!\!0\\
\end{array}\right.\ \ \
\left\{\begin{array}{lll}
y_2C-y_{n+2}C&\!\!\!\!=&\!\!\!\!0\\
y_3C-y_{n+3}C&\!\!\!\!=&\!\!\!\!0\\
&\!\!\!\!\vdots&\\
y_{n+1}C-y_{2n+1}C&\!\!\!\!=&\!\!\!\!0
\end{array}\right.\ \ \
\left\{\begin{array}{lll}
y_1D-y_{n+2}D&\!\!\!\!=&\!\!\!\!0\\
y_{2}D-y_{n+3}D&\!\!\!\!=&\!\!\!\!0\\
&\!\!\!\!\vdots&\\
y_{n-1}D-y_{2n}D&\!\!\!\!=&\!\!\!\!0\\
y_{n}D-y_{2n+1}D&\!\!\!\!=&\!\!\!\!0\\
\end{array}\right.
}\leqno{(\ast\!\ast)^7}
$$

\noindent are equal.

One checks that the system $(\ast\!\ast)^7$ is equivalent to the
matrix equation
$$
[y_{n+1}\ y_{2n+1}\ y_n\ y_{2n}\ y_{n-1}\ y_{2n-1}\ \ldots\ y_2\
y_{n+2}\ y_1]\cdot\CN_n=0,
$$
where $\CN_n=\DM(\ZIdmjj,A,-C,n)$. Therefore the dimension
$[M,I(2n+1,1)]$ equals $\crk(\CN_n)$.

\smallskip

\underline{Case 8: $X=I(2n+1,i)$, $i\in\{2,3,4\}$}. The assertions
for this case follow immediately from the case 7, since the module
$I(2n+1,i)$, for $i\in\{2,3,4\}$, can be obtained from
$I(2n+1,i-1)$ after the cyclic permutation $(1,2,3,4)$ of vertices
of the quiver $Q$ (see Proposition 2).

\smallskip

\underline{Case 9: $X=I(2n, 1)$}. $I(0,1)$ is the injective
$\La$-module
$I(1)=(\trivm_{0,1},\trivm_{0,0},\trivm_{0,0},\trivm_{0,0}\,)$. So
every homomorphism $F:M\to I(0,1)$ is given by the collection
$(0,s,0,0,0)$ where $s\in\MM_{1\times n_1}$ can have an arbitrary
value (see $(*)$). So $[M,I(0,1)]=\dim(\MM_{1\times n_1})=n_1$.

For $n\geq 1$, every homomorphism $F:M\to I(2n,1)$ is given by the
collection $(F_0,F_1,F_2,F_3,F_4)\in \MM_{2n\times
n_0}\times\MM_{(n+1)\times n_1}\times\MM_{n\times
n_2}\times\MM_{n\times n_3}\times\MM_{n\times n_4}$ satisfying the
commutativity relations
$$F_0A=\wmm{\wlj_{n,n+1}}{\wld_{n,n+1}}F_1,\ \ F_0B=\wmm{0_n}{I_n}F_2,\
\ F_0C=\wmm{I_{n}}{0_n}F_3,\ \
F_0D=\wmm{I_{n}}{I_n}F_4\leqno{(\ast)^{9}}$$ (see $(\ast)$ and
Proposition 2).
 Let us denote by $y_1,\ldots,y_{2n}$, $s_1,\ldots, s_{n+1}$,
$t_1,\ldots,t_{n}$, $u_1,\ldots,u_{n}$ and respectively
$w_1,\ldots,w_{n}$, the consecutive rows of matrices $F_0,F_1,$ $
F_2,F_3,F_4$. Similarly as before, the variables $s_i, t_i, u_i,
w_i$ are determined by $y_1,\ldots,$ $y_{2n}$, so clearly the
dimensions of the solution spaces for $(\ast)^9$ and
$${\scriptsize \left\{\begin{array}{lll}
y_{2}A-y_{n+1}A&\!\!\!\!=&\!\!\!\!0\\
y_{3}A-y_{n+2}A&\!\!\!\!=&\!\!\!\!0\\
&\!\!\!\!\vdots&\\
y_{n-1}A-y_{2n-2}A&\!\!\!\!=&\!\!\!\!0\\
y_{n}A-y_{2n-1}A&\!\!\!\!=&\!\!\!\!0\\
\end{array}\right.\  \ \
\left\{\begin{array}{lll}
y_{1}B&\!\!\!\!=&\!\!\!\!0\\
y_{2}B&\!\!\!\!=&\!\!\!\!0\\
&\!\!\!\!\vdots&\\
y_{n}B&\!\!\!\!=&\!\!\!\!0\\
\end{array}\right.\ \ \
\left\{\begin{array}{lll}
y_{n+1}C&\!\!\!\!=&\!\!\!\!0\\
y_{n+2}C&\!\!\!\!=&\!\!\!\!0\\
&\!\!\!\!\vdots&\\
y_{2n}C&\!\!\!\!=&\!\!\!\!0
\end{array}\right.\ \ \
\left\{\begin{array}{lll}
y_1D-y_{n+1}D&\!\!\!\!=&\!\!\!\!0\\
y_{2}D-y_{n+2}D&\!\!\!\!=&\!\!\!\!0\\
&\!\!\!\!\vdots&\\
y_{n-1}D-y_{2n-1}D&\!\!\!\!=&\!\!\!\!0\\
y_{n}D-y_{2n}D&\!\!\!\!=&\!\!\!\!0\\
\end{array}\right.
}\leqno{(\ast\!\ast)^9}
$$

\noindent are equal.

One checks that the system $(\ast\!\ast)^9$ is equivalent to the
matrix equation
$$
[y_{1}\ y_{n+1}\ y_2\ y_{n+2}\ y_{3}\ y_{n+3}\ \ldots\  y_{n}\
y_{2n}]\cdot\CN_n=0,
$$
where $\CN_n=\DM(\ZIdmj,\ZpIdmj,-A,n-1)$. Therefore the dimension
$[M,I(2n,1)]$ equals $\crk(\CN_n)$.

\smallskip

\underline{Case 10: $X=I(2n,i)$, $i\in\{2,3,4\}$}. The assertions
for this case follow immediately from the case 9, since the module
$I(2n,i)$, for $i\in\{2,3,4\}$, can be obtained from $I(2n,i-1)$
after the cyclic permutation $(1,2,3,4)$ of vertices of the quiver
$Q$ (see Proposition 2).

\smallskip

\underline{Case 11: $X=R(l,\la)$}. For $l\geq 1$, every
homomorphism from  $M$ to $R(l,\la)$ is given by the collection of
matrices $(F_0,F_1,F_2,F_3,F_4)\in \MM_{2l\times
n_0}\times\MM_{l\times n_1}\times\MM_{l\times
n_2}\times\MM_{l\times n_3}\times\MM_{l\times n_4}$ satisfying the
commutativity relations $$F_0A=\ma21{I_l}{0_l}F_1,\ \
F_0B=\ma21{0_l}{I_l}F_2,\ \ F_0C=\ma21{I_l}{I_l}F_3,\ \
F_0D=\ma21{J_l(\la)}{I_l}F_4\leqno{(\ast)^{11}}$$ (see $(\ast)$
and Proposition 2). Let us denote by
$y_1,\ldots,y_l,y_{l+1},\ldots,y_{2l}$, $s_1,\ldots, s_l$,
$t_1,\ldots,t_l$, $u_1,\ldots,u_l$ and respectively
$w_1,\ldots,w_l$, the consecutive rows of matrices
$F_0,F_1,F_2,F_3,F_4$. So the system $(\ast)^{11}$ can be
presented in the form
$${\scriptsize \left\{\begin{array}{lll}
y_1A&\!\!\!\!=&\!\!\!\!s_1\\
&\!\!\!\!\vdots&\\
y_lA&\!\!\!\!=&\!\!\!\!s_l\\
y_{l+1}A&\!\!\!\!=&\!\!\!\!0\\
&\!\!\!\!\vdots&\\
y_{2l}A&\!\!\!\!=&\!\!\!\!0\\
\end{array}\right.\  \ \
\left\{\begin{array}{lll}
y_1B&\!\!\!\!=&\!\!\!\!0\\
&\!\!\!\!\vdots&\\
y_lB&\!\!\!\!=&\!\!\!\!0\\
y_{l+1}B&\!\!\!\!=&\!\!\!\!t_1\\
&\!\!\!\!\vdots&\\
y_{2l}B&\!\!\!\!=&\!\!\!\!t_l
\end{array}\right.\ \ \
\left\{\begin{array}{lll}
y_1C&\!\!\!\!=&\!\!\!\!u_1\\
&\!\!\!\!\vdots&\\
y_lC&\!\!\!\!=&\!\!\!\!u_l\\
y_{l+1}C&\!\!\!\!=&\!\!\!\!u_1\\
&\!\!\!\!\vdots&\\
y_{2l}C&\!\!\!\!=&\!\!\!\!u_l
\end{array}\right.\ \ \
\left\{\begin{array}{lll}
y_1D&\!\!\!\!=&\!\!\!\!\la w_1+w_2\\
&\!\!\!\!\vdots&\\
y_{l-1}D&\!\!\!\!=&\!\!\!\!\la w_{l-1}+w_l\\
y_lD&\!\!\!\!=&\!\!\!\!\la w_l\\
y_{l+1}D&\!\!\!\!=&\!\!\!\!w_1\\
&\!\!\!\!\vdots&\\
y_{2l}D&\!\!\!\!=&\!\!\!\!w_l
\end{array}\right.
}\leqno{(\ast\ast)^{11}}
$$
Since again all the variables $s_i, t_i, u_i, w_i$ are determined
by $x_1,\ldots,x_{2l}$, so clearly the dimensions of the solution
spaces for $(\ast\ast)^{11}$ and
$${\scriptsize \left\{\begin{array}{lll}
y_{l+1}A&\!\!\!\!=&\!\!\!\!0\\
y_{l+2}A&\!\!\!\!=&\!\!\!\!0\\
&\!\!\!\!\vdots&\\
y_{2l}A&\!\!\!\!=&\!\!\!\!0\\
\end{array}\right.\  \ \
\left\{\begin{array}{lll}
y_1B&\!\!\!\!=&\!\!\!\!0\\
y_2B&\!\!\!\!=&\!\!\!\!0\\
&\!\!\!\!\vdots&\\
y_lB&\!\!\!\!=&\!\!\!\!0\\
\end{array}\right.\ \ \
\left\{\begin{array}{lll}
y_1C-y_{l+1}C&\!\!\!\!=&\!\!\!\!0\\
y_2C-y_{l+2}C&\!\!\!\!=&\!\!\!\!0\\
&\!\!\!\!\vdots&\\
y_lC-y_{2l}C&\!\!\!\!=&\!\!\!\!0
\end{array}\right.\ \ \
\left\{\begin{array}{lll}
y_1D-\la y_{l+1}D-y_{l+2}D&\!\!\!\!=&\!\!\!\!0\\
&\!\!\!\!\vdots&\\
y_{l-1}D-\la y_{2l-1}D-y_{2l}D&\!\!\!\!=&\!\!\!\!0\\
y_lD-\la y_{2l}D&\!\!\!\!=&\!\!\!\!0\\
\end{array}\right.
}\leqno{(\ast\!\ast\!\ast)^{11}}
$$

\noindent are equal. One easily checks that the system
$(\ast\ast\ast)^{11}$ is equivalent to the matrix equation
$$
[y_l\ y_{2l}\ y_{l-1}\ y_{2l-1}\ y_{l-2}\ y_{2l-2}\ \ldots\  y_1\
y_{l+1}]\cdot\CN_l=0,
$$
where $\CN_l=\DM(\ZRlla,\ZRlla,-D,l-1)$ is the matrix from the
assertion. Therefore the dimension $[M,R(l,\la)]$ is the dimension
of the solution space of the above equation which equals
$\crk(\CN_l)$.

\smallskip

\underline{Case 12: $X=R(s,2l,\la)$, $s\in\ZZ_2$,
$\la\in\{0,1,\infty\}$}. Observe that for every $l\geq 1$,
$R(0,2l,0)$ is just the module $R(l,\la)$ for $\la:=0$ (see
Proposition 2). So we obtain the assertion for $R(0,2l,0)$
immediately from the assertion for $R(l,\la)$ when we substitute
$\la:=0$ (cf. case 11).

The assertions for the remaining modules in this case follow
obviously from the above case of $R(0,2l,0)$ since the modules
$R(s,2l,\la)$ can be obtained from $R(0,2l,0)$ after suitable
permutation of the vertices of quiver $Q$ (see Proposition 2).

\smallskip

\underline{Case 13: $X=R(0,2l-1,0)$}. For $l\geq 1$, every
homomorphism $F:M\to R(0,2l-1,0)$ is given by the collection
$(F_0,F_1,F_2,F_3,F_4)\in \MM_{(2l-1)\times
n_0}\times\MM_{(l-1)\times n_1}\times\MM_{l\times
n_2}\times\MM_{(l-1)\times n_3}\times\MM_{l\times n_4}$ satisfying
the commutativity relations
$$F_0A=\wmmm{I_{l-1}}{0_{1,l-1}}{I_{l-1}}F_1,\ \ F_0B=\wmm{\wlj_{l-1,l}}{I_l}F_2,\
\ F_0C=\wmm{I_{l-1}}{0_{l,l-1}}F_3,\ \
F_0D=\wmm{0_{l-1,l}}{I_l}F_4\leqno{(\ast)^{13}}$$ (see $(\ast)$
and Proposition 2).
 Let us denote by $y_1,\ldots,y_{2l-1}$, $s_1,\ldots, s_{l-1}$,
$t_1,\ldots,t_{l}$, $u_1,\ldots,u_{l-1}$ and respectively
$w_1,\ldots,w_{l}$, the consecutive rows of matrices $F_0,F_1,$ $
F_2,F_3,F_4$. Similarly as before, the variables $s_i, t_i, u_i,
w_i$ are determined by $y_1,\ldots,$ $y_{2l-1}$, so clearly the
dimensions of the solution spaces for $(\ast)^{13}$ and
$${\scriptsize \left\{\begin{array}{lll}
y_{l}A&\!\!\!\!=&\!\!\!\!0\\
y_{1}A-y_{l+1}A&\!\!\!\!=&\!\!\!\!0\\
y_{2}A-y_{l+2}A&\!\!\!\!=&\!\!\!\!0\\
&\!\!\!\!\vdots&\\
y_{l-1}A-y_{2l-1}A&\!\!\!\!=&\!\!\!\!0\\
\end{array}\right.\  \ \
\left\{\begin{array}{lll}
y_{1}B-y_lB&\!\!\!\!=&\!\!\!\!0\\
y_{2}B-y_{l+1}B&\!\!\!\!=&\!\!\!\!0\\
&\!\!\!\!\vdots&\\
y_{l-1}B-y_{2l-2}B&\!\!\!\!=&\!\!\!\!0\\
\end{array}\right.\ \ \
\left\{\begin{array}{lll}
y_{l}C&\!\!\!\!=&\!\!\!\!0\\
y_{l+1}C&\!\!\!\!=&\!\!\!\!0\\
&\!\!\!\!\vdots&\\
y_{2l-1}C&\!\!\!\!=&\!\!\!\!0
\end{array}\right.\ \ \
\left\{\begin{array}{lll}
y_1D&\!\!\!\!=&\!\!\!\!0\\
y_{2}D&\!\!\!\!=&\!\!\!\!0\\
&\!\!\!\!\vdots&\\
y_{l-1}D&\!\!\!\!=&\!\!\!\!0\\
\end{array}\right.
}\leqno{(\ast\!\ast)^{13}}
$$

\noindent are equal.

One checks that the system $(\ast\!\ast)^{13}$ is equivalent to
the matrix equation
$$
[y_{l}\ y_{1}\ y_{l+1}\ y_{2}\ y_{l+2}\ y_{3}\ \ldots\ y_{2l-2}\
y_{l-1}\ y_{2l-1}]\cdot\CN_l=0,
$$
where $\CN_l=\DM(\ZRzdljz,\ZpRzdljz,-B,l-1)$. So the dimension
$[M,R(0,2l,0)]$ equals $\crk(\CN_l)$.

\smallskip

\underline{Case 14: $X=R(s,2l-1,\la)$, $s\in\ZZ_2$,
$\la\in\{0,1,\infty\}$ and $(s,\la)\neq (0,0)$}. The assertions
for this case follow immediately from the case 13 since the
modules $R(s,2l-1,\la)$ can be obtained from $R(0,2l-1,0)$ after
suitable permutation of the vertices of quiver $Q$ (see
Proposition 2).

\smallskip

We considered above all the cases of the assertion so the proof is
finished. \hfill$\Box$

\bigskip

{\bf Remark.} Observe that all the matrices $\CM_*(-,-,-,-)$
appearing in the assertion of Proposition 1 have very specific,
\lq\lq almost block diagonal'' shapes, which is essentially used
in the paper \cite{df}  to decrease the computational complexity
of considered algorithms. We obtained those shapes thanks to
determining the right orders of the equations and the variables
$y_i$ in the matrix equations in the proof above. Determining
those orders was not always an easy task!

\bigskip

\noindent{Faculty
 of Mathematics and Computer Science\\
Nicolaus Copernicus University\\
Chopina 12/18\\
87-100 Toru\' n, Poland}

\noindent {\footnotesize E-mail: amroz@mat.umk.pl}

\end{document}